\documentclass[11pt,a4paper]{article}
\usepackage{amsmath, amsfonts, amssymb,}
\usepackage{a4wide}
\usepackage{parskip}
\usepackage{enumitem}
\usepackage{xcolor}
\usepackage{indentfirst}
\usepackage{pstricks}
\usepackage{latexsym}
\usepackage{amsthm}
\usepackage{cite}

\newcommand{\e}{\epsilon}
\newcommand{\pa} {\partial}
\newcommand{\al} {\alpha}
\newcommand{\ba} {\beta}
\newcommand{\de} {\delta}

\newcommand{\Om} {\Omega}
\newcommand{\De} {\Delta}
\newcommand{\la} {\lambda}

\newcommand{\ds} {\displaystyle}

\newcommand{\RR}{{\mathbb R}}
\newcommand{\te} {\theta}

\newtheorem{thm}{Theorem}[section]
\newtheorem{rem}{Remark}[section]
\newtheorem{lem}{Lemma}[section]
\newtheorem{prop}{Proposition}[section]

\numberwithin{equation}{section}

\begin{document}
		\setlength{\abovedisplayskip}{3pt}
	\setlength{\belowdisplayskip}{3pt}
	\date{}
	{\vspace{0.01in}
		\title{Critical growth fractional Kirchhoff elliptic problems} 
			

		\author{ {\bf Divya Goel$\,^{1,}$\footnote{e-mail: {\tt divya.mat@iitbhu.ac.in}},  Sushmita Rawat $\,^{2,}$\footnote{e-mail: {\tt sushmita.rawat1994@gmail.com}} and 
				K. Sreenadh$\,^{2,}$\footnote{e-mail: {\tt sreenadh@maths.iitd.ac.in}}} \\ $^1\,$Department of Mathematical Sciences, Indian Institute of Technology (BHU),\\ Varanasi 221005, India\\
			$^2\,$Department of Mathematics, Indian Institute of Technology Delhi,\\ Hauz Khas, New Delhi 110016, India.
			}

		\maketitle
		\begin{abstract}
			This article is concerned with the existence and multiplicity of positive weak solutions for the following fractional Kirchhoff-Choquard problem:
			\begin{equation*}
				\begin{array}{cc}
					\displaystyle  M\left( \|u\|^2\right) (-\Delta)^s u = \ds\lambda f(x)|u|^{q-2}u + \left( \int\limits_{\Omega}  \frac{|u(y)|^{2^{*}_{\mu ,s}}}{|x-y|^ \mu}\, dy\right)  |u|^{2^{*}_{\mu ,s}-2}u \;\text{in} \; \Omega,\\
					u  > 0\quad \text{in} \; \Omega, \,\,
					u  = 0\quad \text{in} \; \mathbb{R}^{N}\backslash\Omega,
				\end{array}
			\end{equation*}
			where $\Omega$ is open bounded domain of $\mathbb{R}^{N}$ with $C^2$ boundary, $N > 2s$ and $s \in (0,1)$, here $M$ models Kirchhoff-type coefficient of the form $M(t) = a + bt^{\te-1}$, where $a, b > 0$ are given constants.
			$(-\Delta)^s$ is fractional Laplace operator, $\lambda > 0$ is a real parameter. Using the variational methods, we explore the existence of solution for ${q} \in (1,2^*_s)$ and $\te \geq 1$.
			Here $2^*_s = \frac{2N}{N-2s}$ and $2^{*}_{\mu ,s} = \frac{2N-\mu}{N-2s}$ is the critical exponent in the sense of Hardy-Littlewood-Sobolev inequality.
			\medskip

			\noindent \textbf{Key words:} fractional Laplacian, Hardy-Littlewood-Sobolev critical exponent, Kirchhoff equation, Concentration-compactness.
			
			\medskip
			
			\noindent \textit{2020 Mathematics Subject Classification: 35A15, 35J60, 35J20.}

		\end{abstract}
		\newpage
		\section{Introduction}
		In this work, we study a class of Kirchhoff-type equation for fractional Laplacian with Choquard term. We consider the following problem
	\begin{equation*}    (P_\lambda)\; \left\{  \begin{array}{cc}
			\displaystyle \left( a + b\|u\| ^{2\theta - 2}\right) (-\Delta)^s u = \lambda f(x)|u|^{q-2}u + \left( \int\limits_{\Omega} \frac{|u(y)|^{2^{*}_{\mu ,s}}}{|x-y|^ \mu}\, dy\right)  |u|^{2^{*}_{\mu ,s}-2}u\; \text{in} \; \Omega,\\
			u>0\; \text{in}\; \Omega,\;  u = 0\; \text{in} \; \mathbb{R}^{N}\backslash\Omega,
		\end{array} \right.
	\end{equation*}
	where $M(t)= a+bt^{\theta-1}$, $\Omega$ is open bounded domain of $\mathbb{R}^{N}$ having smooth  boundary, $N > 2s$ with $s\in (0,1)$, $a, b, \theta, \la$ are positive parameters, where $2^{*}_{\mu ,s} = \frac{2N-\mu}{N-2s}$, and  $\|u\|^2=\int\limits_{\mathbb{R}^{2N}} \frac{|u(x)-u(y)|^2}{|x-y|^{N+2s}}\,dxdy$. The function $f(x)$ is a continuous real valued sign changing function.
 Here, $(-\Delta)^s$ is fractional Laplace operator defined as,
	\begin{equation*}
		(-\Delta)^su(x) := C(N, s) \lim \limits_{\e \to 0} \int\limits_{\mathbb{R}^N\backslash B_{\e}(x)}\dfrac{u(x)-u(y)}{|x-y|^{N+2s}}\,dy, \quad x\in \mathbb{R}^N,
	\end{equation*}
	where $C(N,s)$ is the normalization constant and $B_{\e}(x)$ denotes the open ball of radius $\e$ centered at $x$.\\
	Among the nonlocal elliptic equations, Choquard equations have gained attention in mathematics and physics, due to their numerous applications. One of the first applications was given by Pekar in the framework of quantum theory \cite{pekar} and Lieb \cite{choqlieb} used it in the approximation of Hartree-Fock theory. Consider the following equation 
	\begin{equation}\label{th7}
		-\Delta u + V(x) u = (|x|^{-\mu}* F(u))f(u) \text{ in } \RR^N
	\end{equation}
	where $f \in C(\mathbb{R}, \mathbb{R})$ satisfies some growth condition,   $F$ is  the anti-derivative of $f$ and $V$ is the vanishing potential. These type of equations have been studied by Moroz and Schaftingen \cite{Moroz3}. For detailed state of the art research, readers can refer \cite{Moroz4, M.yang} and references therein.\\
	%
	To make this paper more comprehensible for the readers, we provide a brief general introduction to  existence and multiplicity results for   equations involving nonlinear perturbation and critical non-linearity. 	\\
	In nineties, Azorero and Alonso \cite{azorero} studied the  following problem 
	\begin{equation}\label{Ceq1.9}
		(-\De)^s u = \la |u|^{q-2} u + |u|^{2^*_s-1} \text{ in } \Om, \quad  u =0 \text{ in } \RR^N\setminus \Om
	\end{equation}
	for  $s=1,$ and $2<q<2^*= \frac{2N}{N-2}$. They proved the existence of a non-trivial solution for large $\la$.  While for  $s\in (0,1)$, Barrios, Colorado, Servadei and Soria \cite{barrios2} studied the problem \eqref{Ceq1.9}  for $1 < q < 2^*_s $.  For the convex power case $2 < q< 2^*_s$ the existence of the solution is proved using  Mountain-pass theorem, for suitable value of $\la$ depending on dimension $N$.  On the other hand, for the concave case, authors established the  multiplicity of solutions  for small $\la$. These type of existence results for the problems involving the critical Choquard nonlinearity have been established by  Gao and Yang \cite{M.Yang0}. Subsequently, for  similar type of results in whole space, we cite\cite{miyagaki, shang}  and references therein. 
	
	

	Kirchhoff \cite{ kirchhoff} coined the  following ``Kirchhoff" equations
	\begin{align}\label{er2}
		\rho \frac{\pa^2 u}{\pa t^2} -\left(\frac{P_0}{h}+\frac{E}{2L} \int_0^L \bigg| \frac{\pa u}{\pa x}  \bigg|^2~dx \right)\frac{\pa^2 u}{\pa x^2}=0. 
	\end{align}
	This equation stems from  the classical d'Alembert's wave equation by considering the effects of the changes in the length of the strings during vibrations. The parameters in \eqref{er2} have the following meanings: $L$ is the length of the string, $h$ is the area of cross-section, $E$ is the Young modulus of the material, $\rho$ is the mass density
	and $P_0$ is the initial tension.
	
	Fiscella and Valdinoci \cite{fiscella} explored for the first time, stationary Kirchhoff  equation, in bounded smooth domain of $\mathbb{R}^N$, which takes into account the nonlocal aspect of the tension arising 	from nonlocal measurements of the  fractional length of the string. Indeed, they studied the following 
	critical exponent problem 
	\begin{equation} \label{Ceq10}
		\ds \left(  a +b\|u\|^2  \right) (-\Delta)^s u = \la f(x,u) +   |u|^{2^*_s-2}u \text{ in } \Om,\;\; u =0 \text{ in } \RR^N\setminus \Om. 
	\end{equation}
	Here $f$ is a continuous function of sub-critical growth. Authors established the  existence of non trivial solution for large $\la$.    
	Considering the importance of existence of solutions in critical point theory, it is natural to seek analogue of above mentioned results for the Kirchhoff operator. 
	In the last decade, for various non-linearities $f(x, u)$, Kirchhoff problems have been explored by numerous mathematicians.  For instance, Naimen \cite{naimen} studied \eqref{Ceq10} for $s=1$, $f(x,u)= u$,  $a=1$ and $N=3$ and proved the  existence, nonexistence and uniqueness of positive solutions.  Subsequently in \cite{naimen1}, for $f(x,u)= u^q$, Naimen proved the existence of solution when $q \in (2, 4]$ by variational methods, while for $ q \in (4, 6)$, the existence results were obtained by a  cut-off technique. 
	
	On the other hand, 	Li and Liao \cite{Li} studied the  following problem on whole domain 
		\begin{equation} \label{Ceq1.3}
		\ds \left( a + b\|u\| ^{2}\right) (-\Delta)^s u = \la k(x)|u|^{q-2}u +   \mu|u|^{2^*_s-2}u \text{ in } \RR^N
	\end{equation}
	for $s =1$, $a, b > 0$. Here authors registered the existence of  two positive solutions when $2< q < 2^*$, using minimization argument and Mountain-pass theorem, for some $\mu \in (0, \mu^*)$ and for $\la$ large enough, where $k \in L^{\frac{2^*}{2^*-q}}$ and $k \geq 0$. There is a large volume of literature that examines diverse aspects of existence of solution to problems involving Kirchhoff operator. For a sample, we refer reader to some recent articles  \cite{alves, pucci1,pucci2, pucci3, pucci5,  xiang,xiang1, yao, miyagaki1, pucci_xiang, xiang2}.

	 Recently, Goel and Sreenadh \cite{Div} studied the  problem $(P_\la)$ for $s=1,~ 1<q\leq 2$. To prove the multiplicity of solutions for $1< q <2$, authors used the Nehari Manifold technique whereas for $q=2$, existence of solution is obtained using the Mountain-pass theorem.  While for whole domain,
	  Wang and  Xiang \cite{wang} studied the following Kirchhoff-Choquard equation 
	\begin{equation}\label{IM1}   
		\displaystyle \left( a + b\|u\| ^{2\theta-2}\right) (-\Delta)^s u = \lambda f(x)|u|^{q-2}u +\ba  \left( \int\limits_{\RR^N} \frac{|u(y)|^{2^{*}_{\mu ,s}}}{|x-y|^ \mu}\, dy\right) |u|^{2^{*}_{\mu ,s}-2}u \;\;\text{in} \;\RR^N
	\end{equation}
	where $4s \leq \mu < N$, $\theta=2$ and $2< q < 2^*_s$.  For $f \geq 0$,  $a > 0$ and $b$ sufficiently small, they established the existence of  two non-trivial solutions using minimization argument and Mountain-pass theorem for $\la=\ba$, and  $\la $ large enough. Later Liang, Pucci and Zhang \cite{pucci4} proved the same results for $s=1$ case. {  Wang, Hu and Xiang \cite{sing}, investigated the following equation \begin{equation}\label{IM2}   
			\displaystyle \left( a + b\|u\| ^{p\theta-p}\right) (-\Delta)^s_p u = \lambda \frac{f(x)}{u^\ba} + \left( \int\limits_{\RR^N} \frac{g(x)|u|^q}{|x-y|^ \mu}\, dy\right) g(x)u^{q-1}, u>0 \;\;\text{in} \;\RR^N
		\end{equation}
		where $N\geq 2, 1<p<N/s$ with $s \in (0,1), ~1<q<p^*_{\mu,s}$, and $\theta \in [1,2q),~ p^*_{\mu,s}= \frac{p(2N-\mu)}{2(N-ps)}. $ Here authors proved the existence of two non-negative solutions using Nehari manifold approach. In this que,  Rawat and Sreenadh \cite{rawat},   established the multiplicity and regularity of solutions  for \eqref{IM2} in bounded domain  with $ q= 2^*_{\mu,s}$  and  $f= g=1$. They studied the degenerate Kirchhoff problem \eqref{IM2}, using the minimization argument and by approximating the perturbed problem. }


To the best of our knowledge, there is no attempt  till now to check the   existence and multiplicity of solutions to Kirchhoff-Choquard equation for   $1<q<2^*_s,~ \te \geq 1 $. Our aim here is a modest attempt to bridge this gap. The central idea is to  illustrate a unified approach.  Here in this article, we studied the problem $(P_\la)$ and illustrated the existence and multiplicity of solutions.  To be precise, the aforementioned results in references \cite{wang, pucci4} for $\mu >4s$ encouraged us to ask if analogous results  for a Kirchhoff-Choquard equation  for any   $\mu \in (0,N)$  exists. Section 3 and Section 4 seek to show that this is indeed the case. For instance, we prove the following results   in case of $1 \leq \te < 2^*_{\mu,s}$. 
\begin{thm}\label{Cthm1.4}
	Let   $2 < q < 2\te$, then there exist $\Lambda^*> 0$ such that for $N > 4s$, $(P_\la)$ has at least one positive solution, for all $\la \in (0, \Lambda^*)$.
\end{thm}

\begin{thm}\label{Cthm1.3}
	Let   $2\te \leq q < 2^*_s$, then there exists $\Lambda_* > 0$ such that $(P_\la)$ has at least one positive solution for all $\la \geq \Lambda_*$.
\end{thm}
Observe that we get contrasting results for  $2 < q < 2\te$ and $2\te \leq q < 2^*_s$. While for  $\te \geq 2^*_{\mu,s}$, we prove the existence of two positive solutions.  The novelty of this result is that it proves the multiplicity of solutions for any $\mu \in (0,N)$ which is an open problem in \cite{pucci4}. 
	\begin{thm}\label{Cthm1.5}
		Let $2< q < 2^*_s$, there exists $\Lambda_{**} > 0$ such that for $\la > \Lambda_{**}$ 
		\begin{itemize}
			\item [(i)] when $\te = 2^*_{\mu,s}$, $a > 0$ and $b > (S_s^H)^{-2^*_{\mu,s}}$, $(P_\la)$ admits at least two positive solutions,
			\item [(ii)] when $\te > 2^*_{\mu,s}$, there exist $\mathfrak{A}$, $ \mathfrak{B}$ such that either for $a > 0$ and $b > \mathfrak{B}$ or $b >0$ and $a >\mathfrak{A}$, $(P_\la)$ admits at least two positive solutions
		\end{itemize}
	where $\mathfrak{A}:=\frac{\te-2^*_{\mu,s}}{\te-1}\left[\frac{2^*_{\mu,s}-1}{b(\te-1)} \right]^\frac{2^*_{\mu,s}-1}{\te-2^*_{\mu,s}} (S_s^H)^{\frac{-2^*_{\mu,s}(\te-1)}{\te-2^*_{\mu,s}}}$ and $\mathfrak{B}: = \frac{2^*_{\mu,s}-1}{\te-1}\left[\frac{\te-2^*_{\mu,s}}{a(\te-1)} \right]^\frac{\te-2^*_{\mu,s}}{2^*_{\mu,s}-1} (S_s^H)^{\frac{-2^*_{\mu,s}(\te-1)}{2^*_{\mu,s}-1}}$.
	\end{thm}
After this, we deal with the convex-concave behavior of non linearities, i.e. $ 1<q\leq 2$. We have extended the results of \cite{Div} for the fractional diffusion Kirchhoff problems with $\te \geq 1$.  Although this is a fundamental and natural extension, we did not find it  explicitly anywhere in the literature, and for this reason we record it in the last section of the article. To give a consolidated approach, we extend the results of \cite{Div} by using the  minimization approach and Mountain-pass theorem in place of Nehari Manifold technique. Precisely, we prove that 

\begin{thm}\label{Cthm1.1}
	Let $1 \leq \te < 2^*_{\mu,s}$ and $1 < q \leq 2$, then there exist $\Lambda^{**}$, $\tilde{\Lambda}^{**} > 0$ such that 
	\begin{itemize}
		\item[(i)] If $ 0 <\mu < \min\{4s, N\}$, then for $\la \in (0, \Lambda^{**})$ and $1 < q < 2$, $(P_\la)$ admits at least two positive solutions.
		\item  [(ii)]If $ 4s \leq \mu < N$, then for $\la \in (0, \tilde{\Lambda}^{**})$ and $\frac{N}{N-2s}\leq \ q < 2$, $(P_\la)$ admits at least two positive solutions.
		\item[(iii)] For $q= 2$, $(P_\la)$ admits at least one positive solution.
	\end{itemize}
\end{thm}
\begin{rem}\label{Crem1.1}
	Proof of existence of one positive solution in case of $1 < q <2$ holds for all $\te \geq 1$.
\end{rem}
Let us note that the case $M=1$ of $(P_\la)$ are relatively easy to deal as compared to the case $M(t)= a+bt^{\theta-1}$. The presence of the Kirchhoff operator makes the problem more complex to handle as the weak limit of the Palais-Smale sequence is no more a weak solution to the problem. However, our exposition has a different goal: to apply the elementary techniques to establish the above-mentioned results. To demonstrate the proofs, we used a unified approach: the minimization arguments and the Mountain-pass theorem.  
	To the best of our knowledge, this is the first article to address the existence of solution for any $\theta$ and register  contrasting results for $2<q<2\theta$ and $2\theta\leq q <2^*_s$ with $\theta \in [1,2^*_{\mu,s})$.

	\begin{rem}
		The results in this article  can be extended to  problem given on $\RR^N$, like the \cite{pucci4,wang}. 
	\end{rem}
	\begin{rem}
		Also, the conclusions of article can be   generalized to the following  p-fractional problem with critical exponent problem 
		\begin{equation*}    
			\displaystyle \left( a + b\|u\| ^{p\theta - p}\right) (-\Delta)_p^s u = \lambda f(x)|u|^{q-2}u + |u|^{p^{*}_{s}-2}u,\; \text{in} \;\, \Omega,\\
			\;\;u>0\;\; \text{in}\; \Omega,\;  u = 0\quad \text{in} \; \mathbb{R}^{N}\backslash\Omega.
		\end{equation*}
	where $\Om$ is a smooth bounded domain, $\frac{N}{N-ps} < q < p^*_s,~ N>ps,~ a,b,\theta,\la$ are positive parameters, $s \in (0,1]$, $p \geq 2$, $p^*_s= \frac{Np}{N-ps}$, and $f$ is continuous function.  
	\end{rem}
%
%
%
%
%
%
%

Rest of the paper is organized as follows: In section 2, we present  the variational framework of  problem $(P_\la)$.  In section 3, we give some technical lemmas which will help us to prove Theorem \ref{Cthm1.4} and \ref{Cthm1.3} of the paper.  In section 4, we present the proofs of Theorem \ref{Cthm1.4} and \ref{Cthm1.3}.  In section 5, we consider $\te \geq 2^*_{\mu,s}$ and by using Mountain-pass theorem and minimization argument, we obtain two positive solutions for $2< q< 2^*_s$. In section 6, we obtain two positive solutions for the case $1 < q < 2$ and one positive solution for the case $q =2$.

 \section{Preliminaries}
This section targets to set out the background for the current study.
 We define the functional space associated to this problem as
 \begin{equation*}
 	X_0 = \left\lbrace u\in H^s(\mathbb{R}^N): u = 0\; \text{a.e. in} \; \mathbb{R}^N\backslash \Omega  \right\rbrace
 \end{equation*} 
 which is a closed subspace of the fractional Sobolev space $H^s(\mathbb{R}^N)$
 with the corresponding norm,
 \begin{equation*}
 	\|u\|_{X_0} = \|u\| = \displaystyle \left( \int\limits_{\mathbb{R}^{2N}} \frac{|u(x)-u(y)|^2}{|x-y|^{N+2s}}\,dxdy \right)^{\frac{1}{2}}.
 \end{equation*} 
\begin{prop}
	\textbf{(Hardy-Littlewood-Sobolev inequality)}: Let $t$, $r > 1$ and $0 < \mu < N$ with $\frac{1}{t} + \frac{\mu}{N} + \frac{1}{r} = 2$, $f \in L^t(\mathbb{R}^N)$ and $h \in L^r(\mathbb{R}^N)$. Then there exists a sharp constant $C(t, r, \mu, N)$ independent of $f$, $h$ such that
	\begin{equation*}
		\iint\limits_{\mathbb{R}^{2N}} \dfrac{f(x)h(y)}{|x-y|^ \mu}\,dxdy \leq C(t, r, \mu, N)\|f\|_{L^t(\mathbb{R}^N)}\|h\|_{L^r(\mathbb{R}^N)}.
	\end{equation*}
	
\end{prop}

	
From the embedding results \cite{fractS}, we conclude $X_0$ is continuously embedded in $L^{p}(\Omega)$, $ p \in [1, 2_s^{*}]$. Also the embedding is compact for $1 \leq p < 2_s^{*}$. 
The best constant for the embedding $X_0$ into $L^{2_s^{*}}(\mathbb{R}^N)$ is 
\begin{equation}\label{CEQ2.1}
S_s =\inf\limits_{u\in X_0\backslash \{0\}}\left\lbrace\int\limits_{\mathbb{R}^{2N}} \frac{|u(x)-u(y)|^2}{|x-y|^{N+2s}}\,dxdy : \int\limits_{\mathbb{R^N}}|u|^{2_s^{*}} = 1 \right\rbrace.
\end{equation}
Consequently, we define
\begin{equation}\label{Ceq2.2}
S_s^H = \inf\limits_{u\in X_0\backslash \{0\}}\left\lbrace\int\limits_{\mathbb{R}^{2N}} \frac{|u(x)-u(y)|^2}{|x-y|^{N+2s}}\,dxdy:  \int\limits_{\mathbb{R}^{2N}} \dfrac{|u(x)|^{2^{*}_{\mu ,s}}|u(y)|^{2^{*}_{\mu ,s}}}{|x-y|^ \mu}\,dxdy = 1\right\rbrace. 
\end{equation}
We shall summarize briefly the notion and notations of the function where the infimum of  \eqref{CEQ2.1} and \eqref{Ceq2.2} exists but for more details we refer to \cite{Tuhina,fractS}. 
\begin{lem} 
The constant $S_s^H$ is achieved by u if and only if u is of the form\\
$C\left( \frac{t}{t^2 + |x-x_0|^2}\right) ^\frac{N-2s}{2}$, $x \in \mathbb{R}^N$, 
for some $x_0 \in \mathbb{R}^N, C \text{and}\; t > 0.$ Moreover,
$S_s^H = \frac{S_s}{{C(N, \mu)}^\frac{1}{2_{\mu, s}^{*}}}$.
\end{lem}

Consider the family of functions ${U_\epsilon}$, where $U_\epsilon$ is defined as
\begin{equation} \label{minmimizer}
U_\epsilon = \epsilon^{-\frac{N-2s}{2}}u^{*}\left( \frac{x}{\epsilon}\right), \; x\in \mathbb{R}^{N}, \epsilon > 0,
\end{equation}
\begin{equation*}
u^{*}(x) = \overline{u}\left( \frac{x}{{S_s}^\frac{1}{2s}}\right) , \; \overline{u}(x) = \frac{\tilde{u}(x)}{\|\tilde u\|_{L^{2_s^{*}}(\mathbb{R}^N)}} \; \text{and}\; \tilde{u}(x) = \alpha(\beta^2 + |x|^2)^{-\frac{N-2s}{2}},
\end{equation*}
with $\alpha > 0$ and $\beta > 0$ are fixed constants.
Then for each $\epsilon > 0, \; U_\epsilon $ satisfies
\begin{equation*}
(-\Delta)^s u = |u|^{2_s^{*}-2}u \quad in \; \mathbb{R}^N,
\end{equation*}
and the equality,
\begin{equation*}
\int\limits_{\mathbb{R}^{2N}} \frac{|U_\epsilon(x)-U_\epsilon(y)|^2}{|x-y|^{N+2s}}\,dxdy = \int\limits_{\mathbb{R^N}}|U_\epsilon|^{2_s^{*}} = {S_s}^\frac{N}{2s}. 
\end{equation*}

 Without loss of generality, we assume $0 \in \Omega$  and fix $\delta >0$ such that $B_{4\delta} \subset \Omega$. Let $\eta \in C^{\infty}(\mathbb{R}^N)$ be a cut-off function such that
 \begin{equation*}
 	\eta =	\begin{cases}
 		1 & \quad B_{\delta},\\
 		0 & \quad \mathbb{R}^N\backslash B_{2\delta},
 	\end{cases}
 \end{equation*}
 and for each $\epsilon > 0$, let  $u_\epsilon$ be defined as \begin{equation}\label{Ceq2.4}
 	u_\epsilon(x) = \eta(x)U_\epsilon(x)\quad for \; x\in \mathbb{R}^N,
 \end{equation}
 where $U_\epsilon$ is defined in \eqref{minmimizer}.

 \begin{prop}\label{prop2.2}  Let $s \in (0,1)$ and $N > 2s$. Then $$\|u_\epsilon\|^2 \leq S_s^\frac{N}{2s} + O(\epsilon^{N-2s}), \; \|u_\epsilon\|^{2^{*}_s}_{L^{2^{*}_s}} = S_s^{\frac{N}{2s}}+ O(\epsilon^N)\;$$ and
	\begin{equation*}
 		\|u_\epsilon\|_{L^2}^2 \geq \left\{
 		\begin{array}{ll} C_s\epsilon^{2s}+ O(\epsilon^{N-2s}) & N > 4s,\\
 			C_s\epsilon^{2s}|\log(\epsilon)|+ O(\epsilon^{2s}) & N= 4s,\\
 			C_s\epsilon^{N-2s}+ O(\epsilon^{2s}) & N < 4s,\\
 		\end{array} 
 		\right. 
 	\end{equation*}
 	as $\epsilon \to 0$, for some positive constant $C_s$ depending on $s$.
 	
 \end{prop}
 
 \begin{prop}\label{prop2.3} 
 	Let $s \in (0,1)$ and $N > 2s$. Then, the  following  estimates  hold true for some positive constant $C(N, \mu)$
 	\begin{equation*} 
 		\iint\limits_{\Omega\times\Omega} \frac{|u_\e(x)|^{2^{*}_{\mu ,s}}|u_\e(y)|^{2^{*}_{\mu ,s}}}{|x-y|^ \mu}\,dxdy\geq  C(N, \mu)^{\frac{N}{2s}}(S_s^H)^{\frac{2N-\mu}{2s}}-O(\epsilon^N).
 	\end{equation*}
 \end{prop}
Taking into account that we are looking for positive solutions, 
the energy functional associated with the problem $(P_\la)$ is ${J_{\lambda } : X_0(\Omega) \rightarrow \mathbb{R}}$\; defined as,
\begin{equation*}
\begin{aligned}
	\ds J_\lambda (u)= \frac{a}{2}\|u\|^2 + \frac{b}{2\theta}\|u\|^{2\theta} - \frac{\lambda}{q} \int\limits_{\Omega} f(x)(u^+(x))^{q} \,dx 
	-\frac{1}{2\cdot{2^*_{\mu ,s} }}\iint\limits_{\Omega\times\Omega} \frac{(u^+(y))^{2^{*}_{\mu ,s}}(u^+(x))^{2^{*}_{\mu ,s}}}{|x-y|^ \mu}\,dxdy.
\end{aligned}
\end{equation*}
 Using Hardy-Littlewood-Sobolev inequality we can show $J_\la \in C^1$. Indeed for $\phi \in X_0(\Omega)$
\begin{equation}\label{Ceq2.04}
\begin{aligned}
	\displaystyle \langle J'_\lambda (u), \phi\rangle= &\left( a + b \|u\|^{2(\theta-1)}\right)  \int\limits_{\RR^{2N}}\frac{(u(x)-u(y))(\phi(x)-\phi(y))}{|x-y|^{N+2s}}\,dxdy\\
	& - \lambda \int\limits_{\Omega}f(x)(u^+(x))^{q-1}\phi(x)dx
	- \iint\limits_{\Omega\times\Omega}\frac{(u^+(y))^{2^{*}_{\mu ,s}}(u^+(x))^{2^{*}_{\mu ,s}-1}\phi(x)}{|x-y|^ \mu}\,dxdy.
\end{aligned}
\end{equation}
 Throughout the paper we will use the following notations 
 $u^+ = \max\{u, 0\},\, u^- = \max\{-u, 0\} $ and
\begin{equation*}
	\|u\|_{0}^{2\cdot2^{*}_{\mu ,s}} := \iint\limits_{\Omega\times\Omega} \frac{|u(x)|^{2^{*}_{\mu ,s}}|u(y)|^{2^{*}_{\mu ,s}}}{|x-y|^ \mu}\,dxdy.
\end{equation*}

 \section{Technical lemmas}
 We begin by assembling  a couple of vital results,  which will serve as prelude to our main results. Moreover, some results  might  be of independent interest.
 \begin{lem}\label{Clem4.4}
 	Let $u \in X_0(\Om)$ be a non-trivial solution of $(P_\la)$, then $u \in L^\infty(\Om) \cap C^{s}(\RR^N)$. Moreover, $u$ is a positive solution.
 	\proof Let $u $ be a non-trivial solution of $(P_\la)$, assuming $a+b \|u\|^{2\te-2} := \mathcal{C}$(u) where $\mathcal{C}(u) > 0$ since $a,b >0$. Thus the problem $(P_\la)$ can be rewritten as
 	\begin{equation*}
 		\begin{aligned}
 			(-\Delta)^s u &= \ds \frac{\la f(x)(u^+)^{q-1}}{\mathcal{C}(u)}  + \frac{1}{\mathcal{C}(u)} \left( \int\limits_{\Omega} \frac{(u^+(y))^{2^{*}_{\mu ,s}}}{|x-y|^ \mu}\, dy\right)  (u^+)^{2^{*}_{\mu ,s}-1}\; \text{in} \; \Omega,\\
 			u &= 0\quad \text{in} \; \mathbb{R}^{N}\backslash\Omega.
 		\end{aligned}
 	\end{equation*} 
 	{Employing}\cite[Theorem 2.2]{Divya do-reg}, we obtain $u \in L^\infty(\Om) \cap C^{s}(\RR^N)$. Next we note that for $x, y \in\RR^N$, the following holds
 	\begin{equation}\label{Ceq3.20}
 		\begin{aligned}
 			(u(x)- u(y))(u^-(x)-u^-(y)) 
 			&\leq -{|u^-(x)-u^-(y)|}^2.
 		\end{aligned}
 	\end{equation}
 	Thus, by taking $\phi= u^-$ in \eqref{Ceq2.04} and using \eqref{Ceq3.20} we get
 	\begin{equation*}
 		0 = \left( a+b\|u\|^{2\theta -2}\right) \langle u, u^-\rangle \leq -\left( a+b\|u\|^{2\theta -2}\right)\|u^-\|^2.
 	\end{equation*}
 	It implies $u \geq 0$. Therefore, by the maximum-principle, we get $u >0$. \qed
 \end{lem}
The next Lemma shows that the functional satisfies   Mountain-pass geometry for $q \in (1, 2^*_s)$ and $\te \in [1, 2^*_{\mu,s}).$
\begin{lem}\label{lem4.1} 
	The functional $J_\la$ satisfies the following conditions:
	\begin{itemize}
		\item [(i)]There exist $\al,\rho >0$ such that $J_\la(u) > \al$ for $\|u\| =\rho$.
		\item [(ii)]$J_\la(0) = 0$ and there exists $e \in X_0(\Om)$ with $\|e\| > \rho$ and $J_\la(e) < 0$.
	\end{itemize}
	\proof (i). We divide the proof into the following cases:\\
	\textbf{Case 1.} When $2 < q < 2\te$\\
	 Applying H\"older's inequality and Sobolev embedding we get 
	\begin{equation*}
		\begin{aligned}
			J_\la(u) &\geq \frac{a}{2}\|u\|^2 + \frac{b}{2\theta}\|u\|^{2\theta} - \frac{\lambda}{q} \|f\|_{L^r}S_s^{-q/2}\|u\|^{q} 
			-\frac{1}{2\cdot{2^*_{\mu ,s}}}(S_s^H)^{-2^*_{\mu ,s}}\|u\|^{2\cdot{2^*_{\mu ,s}}}\\
			& > \|u\|^{2}\left( \frac{a}{2} -
			\frac{\lambda}{q}\|f\|_{L^r}S_s^{-q/2}\|u\|^{q-2}
			-\frac{1}{2\cdot{2^*_{\mu ,s}}}(S_s^H)^{-2^*_{\mu ,s}}\|u\|^{2\cdot{2^*_{\mu ,s}}-2}\right).
		\end{aligned}
	\end{equation*}
    Since $2  < q < 2\cdot{2^*_{\mu ,s}}$, the function
	\begin{equation*}
		\mathcal{H}(t) := \frac{\lambda}{q}\|f\|_{L^r}S_s^{-q/2}t^{q-2}
		+\frac{1}{2\cdot{2^*_{\mu ,s}}}(S_s^H)^{-2^*_{\mu ,s}}t^{2\cdot{2^*_{\mu ,s}}-2}, \quad t >0,
	\end{equation*} is an increasing function. Thus there exist $\al, \rho > 0$ such that $J_\la(u) > \al$ for $\|u\| = \rho$.\\
	\textbf{Case 2.} When $2\te \leq q < 2^*_s$\\
	Using the same analysis as in Case 1, the result will follow for this range of $q$ as well.\\
	\textbf{Case 3.} When $1< q < 2$\\
	In this case we show that there exists $\la_0$ such that for $\la \in (0, \la_0]$, the functional $J_\la$ satisfies the required geometry. By Case 1. we have
	\begin{equation*} J_\la (u) > \|u\|^{q}\left( \frac{a}{2}\|u\|^{2-q}
		-\frac{1}{2\cdot{2^*_{\mu ,s}}}(S_s^H)^{-2^*_{\mu ,s}}\|u\|^{2\cdot{2^*_{\mu ,s}}-q}-\frac{\lambda}{q}\|f\|_{L^r}S_s^{-q/2}\right).
	\end{equation*}
	Since $ q < 2 <2\cdot{2_{\mu ,s}^{*}}$, the function 
	\begin{equation*}
		\mathcal{G}(t) := {\frac{at^{2-q}}{2}}-{\frac{ (S_s^H)^{-2^{*}_{\mu ,s}} }{2\cdot{2^*_{\mu ,s} }}}{t}^{2\cdot{2^*_{\mu ,s} }-q}, \quad t >0,
	\end{equation*}
	attains maxima at $\rho := \left[\frac{a(2-q)2^*_{\mu,s}(S_s^H)^{2^{*}_{\mu ,s}}}{2\cdot2^{*}_{\mu ,s}-q} \right]^\frac{1}{2\cdot2^*_{\mu,s}-2} $, where $\mathcal{G}(\rho) > 0$.\\ We set $\lambda_0 :=\frac{qS_s^{q/2}}{2\|f\|_{L^r}}\mathcal{G}(\rho),$
	it follows that for any $\lambda \leq \lambda_0$, there exists $\al > 0$ such that $J_\la(u) \geq {\rho^{q} } \frac{\mathcal{G}(\rho)}{2}= \al$ for $\|u\| = \rho$.\\ 
	\textbf{Case 4.} When $q =2$\\
	For $\la \in \left(0, \frac{aS_s}{\|f\|_{L^r}}\right)$, it is easy to see that for small enough $\|u\|$, the functional $J_\la$ satisfies the desired geometry.
	\proof (ii). For every $u \in X_0(\Om)$ and $t> 0$, we have
	\begin{equation*}
		J_\la(tu) = \frac{at^2}{2}\|u\|^2 + \frac{bt^{2\te}}{2\theta}\|u\|^{2\theta} - \frac{\lambda t^q}{q} \int\limits_{\Omega} f(x)(u^+(x))^{q} \,dx 
		-\frac{t^{2\cdot{2^*_{\mu ,s}}}}{2\cdot{2^*_{\mu ,s}}}\|u^+\|_0^{2\cdot{2^*_{\mu ,s}}} \to -\infty,
	\end{equation*}
	as $t \to \infty$. Therefore, we can find a sufficiently large $e \in X_0$, with $\|e\| > \rho$, such that $J_\la(e) < 0$. Hence, the proof follows. \qed
\end{lem}

 Concerning the boundedness of the $(PS)_c$ sequence of the functional $J_\la$ when   $ q \in (1, 2^*_s)$  and $\theta \in [1, 2^*_{\mu,s})$, we have the following result.
\begin{lem}\label{lem3.1}
	Let  $\{u_n\}$ be a $(PS)_c$ sequence for the functional $J_\la$, then $\{u_n\}$ is a bounded sequence in $X_0(\Om)$.
\end{lem}
\proof Let  $\{u_n\} \subset X_0(\Om)$ be a $(PS)_c$ sequence. Then we have
\begin{equation*}
	J_\la(u_n) \rightarrow c \;\; \text{and}\;\; |\langle J_\la'(u_n), \phi\rangle| \leq \epsilon_n \|\phi\|\;\; \text{for all} \;\phi \in X_0(\Om),\; \text{where} \;\;\epsilon_n \to 0\;\; \text{as}\;\; n \to \infty.		
\end{equation*}
Let us assume by contradiction that $\|u_n\| \to \infty $ as $n \to \infty$. We divide the proof into the following cases.\\
\textbf{Case (i)} When $2 < q < 2\te$\\
Using the H\"older's inequality and Sobolev embedding we deduce that
\begin{equation}\label{Cq.1}
	\begin{aligned}
		& J_\la(u_n) - \dfrac{1}{2\cdot2^*_{\mu,s}} \langle J_\la'(u_n), u_n\rangle \\
		& \geq a\left[ \frac{1}{2} - \frac{1}{2\cdot2^*_{\mu,s}}\right] \|u_n\|^2 + b\left[ \frac{1}{2\te} - \frac{1}{2\cdot2^*_{\mu,s}}\right] \|u_n\|^{2\te}-\la \left[ \frac{1}{q} - \frac{1}{2\cdot2^*_{\mu,s}}\right] \frac{\|f\|_{L^r}\|u_n\|^q}{S_s^{q/2}}.
	\end{aligned}
\end{equation}
It implies that
\begin{equation*}
	\frac{c(1+ \e_n\|u_n\|)}{\|u_n\|^q} \geq \frac{a\left[ \frac{1}{2} - \frac{1}{2\cdot2^*_{\mu,s}}\right]}{\|u_n\|^{q-2}} + \frac{b\left[ \frac{1}{2\te} - \frac{1}{2\cdot2^*_{\mu,s}}\right]}{\|u_n\|^{q-2\te}} - \frac{\la \left[ \frac{1}{q} - \frac{1}{2\cdot2^*_{\mu,s}}\right]\|f\|_{L^r}}{S_s^{q/2}},
\end{equation*}
which is not possible. Therefore we get that $\{u_n\}$ is a bounded sequence in $X_0(\Om)$.\\
\textbf{Case (ii)} When $2\te < q < 2^*_s $\\
Consider
\begin{equation*}
	\begin{aligned}
		& J_\la(u_n) - \dfrac{1}{q} \langle J_\la'(u_n), u_n\rangle \\
		& = a\left[ \frac{1}{2} - \frac{1}{q}\right] \|u_n\|^2 + b\left[ \frac{1}{2\te} - \frac{1}{q}\right] \|u_n\|^{2\te}-\left[ \frac{1}{2\cdot2^*_{\mu,s}}- \frac{1}{q}\right] \|u_n^+\|_0^{2\cdot2^*_{\mu,s}}.
	\end{aligned}
\end{equation*}
Let us define $v_n := \dfrac{u_n}{\|u_n\|}$, which implies $\{v_n\} $ is a bounded sequence. So we deduce that
\begin{equation*}
	\frac{c(1+ \e_n\|u_n\|)}{\|u_n\|^{2}} \geq a\left[ \frac{1}{2} - \frac{1}{q}\right] + b\left[ \frac{1}{2\te} - \frac{1}{q}\right]\|u_n\|^{2\te-2} -\left[ \frac{1}{2\cdot2^*_{\mu,s}}- \frac{1}{q}\right] \|u_n\|^{2\cdot2^*_{\mu,s}-2}\|v_n^+\|_0^{2\cdot2^*_{\mu,s}},
\end{equation*}
which gives us a contradiction. Consequently, $\{u_n\}$ is a bounded sequence in $X_0(\Om)$.\\
For $q = 2\te$, we follow the same procedure as in Case (ii) to get the desired result.\\
\textbf{Case (iii)} When $1< q <2$\\
From \eqref{Cq.1}, we imply that
\begin{equation*}
	\frac{c(1+ \e_n\|u_n\|)}{\|u_n\|^2} \geq a\left[ \frac{1}{2} - \frac{1}{2\cdot2^*_{\mu,s}}\right] + b\left[ \frac{1}{2\te} - \frac{1}{2\cdot2^*_{\mu,s}}\right]\|u_n\|^{2\te-2} - \la \left[ \frac{1}{q} - \frac{1}{2\cdot2^*_{\mu,s}}\right]\frac{\|f\|_{L^r}}{S_s^{q/2}}\|u_n\|^{q-2}
\end{equation*}
which is not possible. Therefore we get that $\{u_n\}$ is a bounded sequence in $X_0(\Om)$.\\
\textbf{Case (iv)} When $q =2$\\ 
We follow the same argument as in Case (iii) to get the desired result.\qed \\
Let the minimax value
\begin{equation*}
	c_\la := \inf\limits_{h \in \Gamma}\max\limits_{t \in [0,1]}J_\la(h(t)),
\end{equation*}
where 
\begin{equation*}
	\Gamma = \{h \in C([0,1], X_0(\Om)): h(0) = 0\;\text{and}\; h(1) = e\}.
\end{equation*} 
\begin{prop}\label{prop4.1}
	Let $2<q < 2\te$, $\te \in [1,2^*_{\mu,s})$ and $\{u_n\}$ be a $(PS)_c$ for $J_\la$ with 
	\begin{equation*}
		-\infty <c < c^*_1 :=\left( \frac{1}{2\te}- \frac{1}{2\cdot{2^*_{\mu ,s}}}\right)b^{\frac{2^*_{\mu,s}}{2^*_{\mu,s}-\te}} (S_s^H)^{\frac{2^*_{\mu,s}\te}{2^*_{\mu,s}-\te}}-\la^{\frac{2\te}{2\te-q}}\hat{D_1},
	\end{equation*}
	where $\hat{D_1} = \ds\frac{(2\te-q)}{(2\te)}\left(\frac{( 2\cdot2^*_{\mu,s}-q)\|f\|_{L^r} }{2q\cdot2^*_{\mu,s}S_s^{\frac{q}{2}}} \right)^{\frac{2\te}{2\te-q}}\left( \frac{2^*_{\mu,s}\cdot q}{b(2^*_{\mu,s}-\te)}\right)^{\frac{q}{2\te-q}}.$
	Then $\{u_n\}$ contains a convergent subsequence.
	\proof Let $\{u_n\}$ be a $(PS)_c$ sequence for $J_\la$ then by Lemma \ref{lem3.1}, we obtain $\{u_n\}$ is a bounded sequence in $X_0(\Om)$. Thus, there exists $u \in X_0(\Om)$ such that up to a subsequence,
	\begin{equation}\label{Ceq4.2}\left\{
		\begin{array}{lr}
			u_n\rightharpoonup u\; \text{in} \; X_0(\Om),\; \|u_n\|\rightarrow \ba,\;
			u_n\rightarrow u\; \text{in} \; L^{p}(\Omega)\; \text{for all}\; {p\in [1,2_s^*)},\\
			\ds\left( \int\limits_{\Omega} \frac{(u_n^+(y))^{2^{*}_{\mu ,s}}}{|x-y|^ \mu}\,dy\right) (u_n^+)^{2^{*}_{\mu ,s}-1} \rightharpoonup \left( \int\limits_{\Omega} \frac{(u^+(y))^{2^{*}_{\mu ,s}}}{|x-y|^ \mu}\,dy\right) (u^+)^{2^{*}_{\mu ,s}-1}\; \text{weakly in} \; L^\frac{2N}{N+2s},\\
			\|u_n^+-u^+\|_0\rightarrow d ,\\
			u_n\rightarrow u \quad\text{a.e. in}\; \Om.
		\end{array}
		\right.
	\end{equation}
	Using \eqref{Ceq4.2} we deduce that
	\begin{equation}\label{Ceq4.3}
		\begin{aligned}
			o(1) =& \left\langle J'_\la(u_n), u_n-u \right\rangle \\
			=& \left( a + b\|u_n\|^{2\theta - 2}\right) \left\langle u_n, u_n-u \right\rangle -\la\int\limits_{\Om}f(x)(u_n^+(x))^{q-1}(u_n-u)(x)\,dx\\
			& -\iint\limits_{\Om\times\Om}\frac{(u_n^+(y))^{2^*_{\mu,s}}(u_n^+(x))^{2^*_{\mu,s}-1}(u_n-u)(x)}{|x-y|^\mu}\,dxdy.
		\end{aligned}
	\end{equation}
	Employing \eqref{Ceq4.2} and Br\'ezis-Lieb
	lemma \cite[Theorem 1]{BrezLieb}, we get as $n \to \infty$
	\begin{equation}\label{Ceq4.4}
		\left( a + b \|u_n\|^{2\theta - 2}\right) \left\langle u_n, u_n-u \right\rangle = \left( a + b \ba  ^{2\theta - 2}\right)\left( \ba^2-\|u\|^2\right) =\left( a + b \ba  ^{2\theta - 2}\right)\|u_n-u\|^2 + o(1).
	\end{equation}  
	By Riesz representation theorem and \eqref{Ceq4.2} we deduce that
	\begin{equation}\label{Ceq4.5}
		\la\int\limits_{\Om}f(x)(u_n^+(x))^{q-1}(u_n-u)(x)\,dx = o(1)\;\;\text{as}\; n \to \infty.
	\end{equation}
	Using Br\'ezis-Lieb lemma \cite{M.yang} and \eqref{Ceq4.2}, we get
		\begin{equation}\label{Ceq4.6}
		\begin{aligned}
			\iint\limits_{\Om\times\Om}\frac{(u_n^+(y))^{2^*_{\mu,s}}(u_n^+(x))^{2^*_{\mu,s}-1}(u_n-u)(x)}{|x-y|^\mu}\,dxdy &= \|u_n^+\|_0^{2\cdot2^*_{\mu,s}}-\|u^+\|_0^{2\cdot2^*_{\mu,s}} + o(1)\\
			&= \|u_n^+-u^+\|_0^{2\cdot2^*_{\mu,s}} + o(1)\\
			&= d^{2\cdot2^*_{\mu,s}}.
		\end{aligned}
	\end{equation} 
	Thus, putting together \eqref{Ceq4.4}, \eqref{Ceq4.5} and \eqref{Ceq4.6} in \eqref{Ceq4.3}, we get as $n \to \infty$
	\begin{equation} \label{Ceq4.7}
		\left( a + b \ba  ^{2\theta - 2}\right)\|u_n-u\|^2 + o(1) = d^{2\cdot2^*_{\mu,s}}. 
	\end{equation}
	If $d = 0$, then $u_n \to u$ strongly in $X_0(\Om)$. So, let us assume $d > 0$.\\
	Deploying Sobolev embedding and \eqref{Ceq4.7}, we obtain
	\begin{equation}\label{Ceq4.8}
		\left( a + b \ba  ^{2\theta - 2}\right)S_s^H \leq d^{2\cdot2^*_{\mu,s}-2}.
	\end{equation}
	From \eqref{Ceq4.7} and \eqref{Ceq4.8}, we deduce that
	\begin{align}\label{Ceq4.9}
		\left[ \left( a + b \ba  ^{2\theta - 2}\right)S_s^H\right]^\frac{2\cdot2^*_{\mu,s}}{2\cdot2^*_{\mu,s}-2}  &\leq d^{2\cdot2^*_{\mu,s}} = \left( a + b \ba  ^{2\theta - 2}\right)\left( \ba^2- \|u\|^2\right) \nonumber\\
		\left( a + b \ba  ^{2\theta - 2}\right)^\frac{1}{2^*_{\mu,s}-1}\left( S_s^H\right)^\frac{2^*_{\mu,s}}{2^*_{\mu,s}-1} &\leq \ba^2- \|u\|^2 < \ba^2 \nonumber\\
		\left( b \ba  ^{2\theta - 2}\right)^\frac{1}{2^*_{\mu,s}-1}\left( S_s^H\right)^\frac{2^*_{\mu,s}}{2^*_{\mu,s}-1} & < \ba^2 \nonumber\\
		 b^\frac{1}{2^*_{\mu,s}-\te}\left( S_s^H\right)^\frac{2^*_{\mu,s}}{2^*_{\mu,s}-\te} & < \ba^2.
	\end{align}
	Arguing as in \eqref{Ceq4.9}, we conclude from \eqref{Ceq4.7}, \eqref{Ceq4.8} and \eqref{Ceq4.9} 
	\begin{equation}\label{Ceq4.10}
		\begin{aligned}
			\left( a + b \ba^{2\theta - 2}\right)^\frac{1}{2^*_{\mu,s}-1}\left( S_s^H\right)^\frac{2^*_{\mu,s}}{2^*_{\mu,s}-1} &\leq \|u_n-u\|^2 + o(1)\\
			\left( b \ba ^{2\theta - 2}\right)^\frac{1}{2^*_{\mu,s}-1}\left( S_s^H\right)^\frac{2^*_{\mu,s}}{2^*_{\mu,s}-1} & < \|u_n-u\|^2 + o(1)\\
			b^{\frac{1}{2^*_{\mu,s}-\te}}\left( S_s^H\right)^\frac{2^*_{\mu,s}}{2^*_{\mu,s}-\te} & < \|u_n-u\|^2 + o(1).
		\end{aligned} 
	\end{equation}
	Employing H\"older’s inequality, Sobolev embedding and Young’s inequality, we get
	\begin{equation}\label{Ceq4.11}
		\begin{aligned}
			\la & \left(\frac{1}{q} - \frac{1}{2\cdot2^*_{\mu,s}} \right)\int\limits_{\Om}f(x)(u^+(x))^q\,dx \leq \, \la\left(\frac{1}{q} - \frac{1}{2\cdot2^*_{\mu,s}} \right)\|f\|_{L^r}\|u\|^qS_s^{-{q/2}}\\
			&= \left( \left[\frac{\te b}{q}\left(\frac{1}{\te} - \frac{1}{2^*_{\mu,s}} \right)  \right]^\frac{q}{2\te}\|u\|^q \right) \left(\left[\frac{\te b}{q}\left(\frac{1}{\te} - \frac{1}{2^*_{\mu,s}} \right)  \right]^\frac{-q}{2\te}\la\left(\frac{1}{q} - \frac{1}{2\cdot2^*_{\mu,s}} \right)\|f\|_{L^r}S_s^{-q/2} \right)\\
			&\leq \, b\left(\frac{1}{2\te} - \frac{1}{2\cdot2^*_{\mu,s}} \right) \|u\|^{2\te} + \la^\frac{2\te}{2\te-q}\hat{D_1},
		\end{aligned}
	\end{equation}
	where $\hat{D_1} = \ds
	\frac{(2\te-q)}{(2\te)}\left(\frac{( 2\cdot2^*_{\mu,s}-q)\|f\|_{L^r} }{2q\cdot2^*_{\mu,s}S_s^{q/2}} \right)^{\frac{2\te}{2\te-q}}\left( \frac{2^*_{\mu,s}q}{b(2^*_{\mu,s}-\te)}\right)^{\frac{q}{2\te-q}}.$\\
	Further from \eqref{Ceq4.2}, \eqref{Ceq4.10}, \eqref{Ceq4.11} and Br\'ezis-Lieb lemma, we get
	\begin{equation*}
		\begin{aligned}
			c \geq&\, J_\la(u_n) - \frac{1}{2\cdot2^*_{\mu,s}}\left\langle J'_\la(u_n), u_n \right\rangle\\
			\geq&\, a\left(\frac{1}{2} - \frac{1}{2\cdot2^*_{\mu,s}} \right)\|u_n\|^2 + b\left(\frac{1}{2\te} - \frac{1}{2\cdot2^*_{\mu,s}} \right)\|u_n\|^{2\te} - \la\left(\frac{1}{q} - \frac{1}{2\cdot2^*_{\mu,s}} \right)\frac{\|f\|_{L^r}\|u\|^q}{S_s^{q/2}} + o(1)\\
			\geq&\, a\left(\frac{1}{2} - \frac{1}{2\cdot2^*_{\mu,s}} \right)\ba^2 + b\left(\frac{1}{2\te} - \frac{1}{2\cdot2^*_{\mu,s}} \right)\left( \|u_n-u\|^{2\te}+\|u\|^{2\te}\right)\\  
			&- \la\left(\frac{1}{q} - \frac{1}{2\cdot2^*_{\mu,s}} \right)\frac{\|f\|_{L^r}\|u\|^q}{S_s^{q/2}} + o(1)\\
			> &\,
			b^{\frac{2^*_{\mu,s}}{2^*_{\mu,s}-\te}}\left(\frac{1}{2\te} - \frac{1}{2\cdot2^*_{\mu,s}} \right)\left(S_s^H\right)^\frac{2^*_{\mu,s}\te}{2^*_{\mu,s}-\te} - \la^\frac{2\te}{2\te-q}\hat{D_1} = c_1^*
		\end{aligned}
	\end{equation*}
	which is not true. Hence $d = 0$ and from \eqref{Ceq4.7}, $u_n \to u$ strongly in $X_0(\Om)$.\qed
\end{prop} 
 \begin{prop}\label{prop3.1}
 	Let $2\te \leq q < 2^*_s$, $\te \in [1,2^*_{\mu,s})$ and $\{u_n\}$ be a $(PS)_c$ for $J_\la$ with 
 	\begin{equation*}
 		-\infty <c < c^*_2 :=\left( \frac{1}{2}- \frac{1}{2\cdot2^*_{\mu,s}}\right) (aS_s^H)^{\frac{2N-\mu}{N-\mu+2s}}.
 	\end{equation*}
 	Then $\{u_n\}$ contains a convergent subsequence.
 	\proof Let $\{u_n\}$ be a $(PS)_c$ sequence for $J_\la$ then by Lemma \ref{lem3.1}, we have $\{u_n\}$ is a bounded sequence. Thus, there exists $u \in X_0(\Om)$ such that up to a subsequence $u_n \rightharpoonup u$ weakly in $X_0(\Om)$, $u_n \to u$ in $L^p(\Om)$ for $p \in [1, 2^*_s)$, $u_n \to u$ a.e. in $\Om$ and $\|u_n\| \to \al$ as a real sequence. 
 	Furthermore, there exist bounded non-negative Radon measures $\omega, \xi$ and $\nu$ such that as $n \to \infty$
 	\begin{equation}\label{Ceq2.11}
 		|(-\De)^ \frac{s}{2}u_n|^2\rightharpoonup \omega, \;|u_n|^{2^*_s} \rightharpoonup \xi \;\;\text{and}\;\; (I_\al * (u_n^+)^{2^*_{\mu,s}})(u_n^+)^{2^*_{\mu,s}} \rightharpoonup \nu 
 	\end{equation} weakly in the sense of measures.  Hence, by \cite[Lemma 3.1]{frchoq}
 	there exist an at-most countable set $I$, a sequence of distinct points $\{x_i\}_{i\in I} \subset \RR^N$ and a family of positive numbers $\{\nu_i\}_{i\in I}$ such that the following holds.
 	\begin{equation}\label{Ceq2.10}
 		\begin{aligned}
 			\nu &= (I_\al * |u|^{2^*_{\mu,s}})|u|^{2^*_{\mu,s}}+ \sum\limits_{i \in I}\nu_i\de_{x_i}, \;\; \sum\limits_{i \in I}\nu_i^{\frac{N}{N+\al}} < \infty;\\
 			\xi &\geq |u|^{2^*_s} + \sum\limits_{i \in I}\xi_i\de_{x_i},\;\; \xi_i \geq C_N (\al)^\frac{-N}{N+\al}\nu_{i}^\frac{N}{N+\al}\;\; \text{and}\\
 			\omega &\geq |(-\De)^\frac{s}{2}u|^2+ \sum\limits_{i \in I}\omega_i\de_{x_i},\;\; \omega_i \geq S_s^H\nu_i^{\frac{1}{2^*_{\mu,s}}},
 		\end{aligned}
 	\end{equation}
 	where $\de_x$ is the Dirac-mass of mass 1 concentrated at $x \in \RR^N$.
 	Let $\e> 0$, we fix a smooth cut-off function $\phi_{\e,i}$ centred at $x_i$ such that 
 	\begin{equation*}
 		0 \leq \phi_{\e,i} \leq 1,\; \phi_{\e,i} \equiv 1 \;\text{in}\; B(x_i, \e/2)\; \text{and}\; \phi_{\e,i} \equiv 0 \;\text{in}\; \RR^N\backslash B(x_i, \e).
 	\end{equation*}
 	Then, by dominated convergence theorem we have
 	\begin{equation*}
 		\int\limits_{\Om}f(x)(u_n^+(x))^q\phi_{\e,i}(x)\;dx  \to \int\limits_{\Om}f(x)(u^+(x))^q\phi_{\e,i}(x)\;dx, \;\text{as}\; n \to \infty
 	\end{equation*}
 	and as $\e \to 0$ we have
 	\begin{equation}\label{Ceq2.12}
 		\lim\limits_{\e \to 0}\lim\limits_{n \to \infty}\int\limits_{\Omega}f(x)(u_n^+(x))^q\phi_{\e,i}(x)\;dx = 0.
 	\end{equation}
 	Consider
 	\begin{align}\label{Ceq2.13}
 		& \left\langle J'_\la(u_n), \phi_{\e,i}u_n \right\rangle \nonumber \\
 		= & \left(a+ b\|u_n\|^{2\te-2} \right) \int\limits_{\RR^{2N}}\frac{\phi_{\e,i}(x)|u_n(x)-u_n(y)|^2}{|x-y|^{N+2s}} + \frac{u_n(y)(u_n(x)-u_n(y))(\phi_{\e,i}(x)- \phi_{\e,i}(y))}{|x-y|^{N+2s}}\\   
 		&- \la\int\limits_{\Omega}f(x)(u_n^+(x))^q\phi_{\e,i}(x)\;dx - \int_{\Om\times\Om}\frac{(u_n^+(y))^{2^*_{\mu,s}}(u_n^+(x))^{2^*_{\mu,s}}\phi_{\e,i}(x)}{|x-y|^\mu}\,dxdy \nonumber.
 	\end{align}
 	Using \eqref{Ceq2.10} and the weak convergence of measure  \eqref{Ceq2.11}, we deduce
 	\begin{align}\label{Ceq2.14}
 		\lim\limits_{\e \to 0}\lim\limits_{n \to \infty}\left(a+ b\|u_n\|^{2\te-2} \right)\int\limits_{\RR^{2N}}\frac{\phi_{\e,i}(x)|u_n(x)-u_n(y)|^2}{|x-y|^{N+2s}} &= \lim\limits_{\e \to 0}\left(a+ b\al^{2\te-2} \right)\int\limits_{\RR^{2N}}\phi_{\e,i}(x)d\omega \nonumber\\
 		&\geq a \omega_i.
 	\end{align}
 	Employing H\"older's inequality, we get
 	\begin{equation*}
 		\left| \int\limits_{\RR^{2N}}\frac{u_n(y)(u_n(x)-u_n(y))(\phi_{\e,i}(x)- \phi_{\e,i}(y))}{|x-y|^{N+2s}}\right| \leq \|u_n\| \left( \int\limits_{\RR^{2N}}\frac{|u_n(y)|^2|\phi_{\e,i}(x)- \phi_{\e,i}(y)|^2}{|x-y|^{N+2s}}\right)^\frac{1}{2}.
 	\end{equation*}
 	As in \cite[Lemma 2.1]{amb} we have
 	\begin{equation}\label{Ceq2.15}
 		\lim\limits_{\e \to 0}\limsup\limits_{n \to \infty}\int\limits_{\RR^{2N}}\frac{|u_n(y)|^2|\phi_{\e,i}(x)- \phi_{\e,i}(y)|^2}{|x-y|^{N+2s}} = 0.
 	\end{equation}
 	Again by \eqref{Ceq2.10} and \eqref{Ceq2.11} we deduce
 	\begin{align}\label{Ceq2.16}
 		\lim\limits_{\e \to 0}\lim\limits_{n \to \infty}\int_{\Om\times\Om}\frac{(u_n^+(y))^{2^*_{\mu,s}}(u_n^+(x))^{2^*_{\mu,s}}\phi_{\e,i}(x)}{|x-y|^\mu}dxdy &= \lim\limits_{\e \to 0}\int\limits_{\Om}\phi_{\e,i}(x)d\nu \nonumber\\
 		&= \nu_i.
 	\end{align}
 	Taking into account \eqref{Ceq2.14}, \eqref{Ceq2.15}, \eqref{Ceq2.12} and \eqref{Ceq2.16} in \eqref{Ceq2.13}, we obtain
 	\begin{equation*}
 		0 = \lim\limits_{\e \to 0}\lim\limits_{n \to \infty}\left\langle J'_\la(u_n), \phi_{\e,i}u_n \right\rangle
 		\geq a\omega_i  - \nu_{i}.
 	\end{equation*}
 	It implies $a\omega_i \leq \nu_{i}$. Combining this with the fact that $ S_s^H \nu_{i}^\frac{1}{2^*_{\mu,s}} \leq \omega_i$ we obtain
 	\begin{equation}\label{Ceq2.17}
 		\omega_i \geq \left( a (S_s^H)^{2^*_{\mu,s}} \right)^\frac{1}{2^*_{\mu,s}-1} \,\,\text{or}\;\; \omega_i = 0.
 	\end{equation}
 	If possible, let there exists $i_0 \in I$ such that $ \omega_{i_0} \geq \left( a (S_s^H)^{2^*_{\mu,s}} \right)^\frac{1}{2^*_{\mu,s}-1} $. Then using the fact that $\{u_n\}$ is a $(PS)_c$ sequence, we deduce
 	\begin{equation*}
 		\begin{aligned}
 			c &= \lim\limits_{n \to \infty}J_\la(u_n)-\frac{1}{q}\left\langle J'_\la(u_n), u_n \right\rangle \hspace{10cm}\\
 			&\geq \, a\left( \frac{1}{2}- \frac{1}{q}\right) \left( \|u\|^2+ \sum\limits_{i \in I}\omega_i\right) + b\left( \frac{1}{2\te}- \frac{1}{q}\right) \left( \|u\|^2+ \sum\limits_{i \in I}\omega_i\right)^\te \\
 			& \quad+ \left( \frac{1}{q}- \frac{1}{2\cdot2^*_{\mu,s}}\right)\left( \|u^+\|_0^{2\cdot2^*_{\mu,s}}+ \sum\limits_{i \in I}\nu_i\right) \\
 			&\geq  \left( \frac{1}{2}- \frac{1}{2\cdot2^*_{\mu,s}}\right) a\omega_{i_0}
 			\geq  \left( \frac{1}{2}- \frac{1}{2\cdot2^*_{\mu,s}}\right) (aS_s^H)^{\frac{2N-\mu}{N-\mu+2s}} = c_2^*\, > c.
 		\end{aligned}
 	\end{equation*}
 	Therefore, $\omega_i = 0$ for all $i \in I$. Hence, we get
 	$$\|u_n\|_0^{2\cdot2^*_{\mu,s}} \to \|u\|_0^{2\cdot2^*_{\mu,s}}.$$
 	Taking into account $J'_\la(u_n) \to 0$ and Br\'ezis-Lieb
 	lemma \cite[Theorem 1]{BrezLieb} we have
 	\begin{equation*}
 		o(1) =  \left\langle J'_\la(u_n), u_n-u \right\rangle = \left(a+ b\al^{2\te-2} \right)\|u_n-u\|^2.
 	\end{equation*}
 	Hence $u_n \to u $ in $X_0(\Om)$, finishing the proof.\qed
 	
 \end{prop}
 
 \begin{prop}\label{P.S}
 	Let $1 < q < 2$, $\te \in [1,2^*_{\mu,s})$ and $\{u_n\}$ be a $(PS)_c$ sequence for $J_\la$ with
 	\begin{equation*}
 		-\infty <c < c^*_3 :=\frac{N-\mu+2s}{2(2N-\mu)}(aS_s^H)^{\frac{2N-\mu}{N-\mu+2s}} - \hat D\la^{\frac{2}{2-q}},
 	\end{equation*}
 	where $\hat D = \frac{(2-q)(2\te-q)}{4\te q}\left( \frac{2\te -q}{2aS_s(\te-1)}\right)^{\frac{q}{2-q}}\|f\|_{L^r}^\frac{2}{2-q}. $ Then $\{u_n\}$ contains a convergent subsequence.
 	\proof Let $\{u_n\}$ be a $(PS)_c$ sequence for $J_\la$ then by Lemma \ref{lem3.1}, $\{u_n\}$ is a bounded sequence. 
 	Using the same arguments up to \eqref{Ceq2.17} as in the proof of Proposition \ref{prop3.1}, we have 
 	\begin{equation}\label{Ceq6.05}
 		\omega_i \geq \left( a (S_s^H)^{2^*_{\mu,s}} \right)^\frac{1}{2^*_{\mu,s}-1} \,\,\text{or}\;\; \omega_i = 0.
 	\end{equation}
 	If possible, let there exists $i_0 \in I$ such that $ \omega_{i_0} \geq \left( a (S_s^H)^{2^*_{\mu,s}} \right)^\frac{1}{2^*_{\mu,s}-1} $.\\
 	Taking into account H\"older’s inequality, Sobolev embedding and Young’s inequality, we get
 	\begin{equation}\label{Ceq2.18}
 		\begin{aligned}
 			\la\int\limits_{\Om}f(x)(u^+(x))^q\,dx &\leq \la \|f\|_{L^r}\|u\|^q (S_s)^{-q/2}\\
 			&\leq \frac{a(\te-1)}{2\te}\left[\frac{1}{q}-\frac{1}{2\te} \right]^{-1}\|u\|^2 + \la^{\frac{2}{2-q}}\frac{2-q}{2}\left[\frac{2\te-q}{2aS_s(\te-1)} \right]^{\frac{q}{2-q}}\|f\|_{L^r}^\frac{2}{2-q}.
 		\end{aligned}
 	\end{equation}
 	Using \eqref{Ceq6.05} and \eqref{Ceq2.18} we deduce that
 	\begin{equation*}
 		\begin{aligned}
 			c =& \lim\limits_{n \to \infty}J_\la(u_n)-\frac{1}{2\te}\left\langle J'_\la(u_n), u_n \right\rangle\\
 			\geq & \left( \frac{1}{2}- \frac{1}{2\te}\right) a\left( \|u\|^2+ \sum\limits_{i \in I}\omega_i\right) - \la \left( \frac{1}{q}- \frac{1}{2\te}\right)\int\limits_{\Omega}f(x)(u^+(x))^{q}dx \\
 			& + \left( \frac{1}{2\theta}- \frac{1}{2\cdot2^*_{\mu,s}}\right)\left( \|u^+\|_0^{2\cdot2^*_{\mu,s}}+ \sum\limits_{i \in I}\nu_i\right) \\
 			\geq &\; \frac{N-\mu+2s}{2(2N-\mu)}(aS_s^H)^{\frac{2N-\mu}{N-\mu+2s}} - \hat D\la^{\frac{2}{2-q}} \,= c_3^*> c.
 		\end{aligned}
 	\end{equation*}
 	
 	Thus, concluding as in Proposition \ref{prop3.1}, we get $u_n \to u $ in $X_0(\Om)$. 
 	\qed

 \end{prop}

\begin{prop}\label{prop3.3}
	Let us assume $\la \in \left(0, \frac{aS_s}{\|f\|_{L^r}}\right)$, $q = 2$, $\te \in [1,2^*_{\mu,s})$ and $\{u_n\}$ be a $(PS)_c$ for $J_\la$ with 
	\begin{equation*}
		-\infty <c < c_4^* := \left( \frac{1}{2}- \frac{1}{2\cdot2^*_{\mu,s}}\right) (aS_s^H)^{\frac{2N-\mu}{N-\mu+2s}}.
	\end{equation*}
	Then $\{u_n\}$ contains a convergent subsequence.
	\proof Let $\{u_n\}$ be a $(PS)_c$ sequence for $J_\la$ then following the same proof as in Proposition \ref{prop3.1}, upto \eqref{Ceq2.17} we have $\omega_i \geq \left( a (S_s^H)^{2^*_{\mu,s}} \right)^\frac{1}{2^*_{\mu,s}-1} \,\,\text{or}\;\; \omega_i = 0.$
	Let us assume, there exists $i_0 \in I$ such that $ \omega_{i_0} \geq \left( a (S_s^H)^{2^*_{\mu,s}} \right)^\frac{1}{2^*_{\mu,s}-1} $.
	Furthermore by Palais-Smale condition we have,
	\begin{equation*}
		\begin{aligned}
			c =& \lim\limits_{n \to \infty}J_\la(u_n)-\frac{1}{2\te}\left\langle J'_\la(u_n), u_n \right\rangle\\
			\geq &\, a\left( \frac{1}{2}- \frac{1}{2\te}\right)  \sum\limits_{i \in I}\omega_i+\left( \frac{1}{2}- \frac{1}{2\te}\right)\left( a -\la\|f\|_{L^r}S_s^{-1}\right) \|u\|^2\\
			& + \left( \frac{1}{2\te}- \frac{1}{2\cdot2^*_{\mu,s}}\right)\left( \|u^+\|_0^{2\cdot2^*_{\mu,s}}+ \sum\limits_{i \in I}\nu_i\right) \\
			\geq &\; \left( \frac{1}{2}- \frac{1}{2\cdot2^*_{\mu,s}}\right) (aS_s^H)^{\frac{2N-\mu}{N-\mu+2s}}\, = c_4^*> c.
		\end{aligned}
	\end{equation*}
	Therefore, $\omega_i = 0$ for all $i \in I$. Hence, 
	concluding as in Proposition \ref{prop3.1} we get our desired result.
\end{prop}
 \section{Proof  of Theorem \ref{Cthm1.4} and \ref{Cthm1.3}}
This section is devoted to the proof of the
	Theorem \ref{Cthm1.4} and Theorem \ref{Cthm1.3}. At the outset, we give some standard observations.\\
  Since, $f$ is a continuous function on $\Om$ and $f^+ = \max\{f(x), 0\}\neq 0$, the set $\Xi
 = \{x \in \Om: f(x) > 0\}$ is an open set of positive measures. Without loss of generality, let us assume $\Xi$ is a domain and $0 \in \Xi$. This implies there exists a $\de > 0$ such that $B_{4\de}(0) \subset \Xi \subseteq \Om$ and $f(x) > 0$ for all $x \in B_{2\de}(0)$. It implies that there exists a $m_f > 0$ such that $f(x) > m_f$ for all $x \in B_{2\de}(0)$. Next choose $\la_1 > 0$ such that for $\la \in (0, \la_1)$, $c_1^* > 0$ (as defined in Proposition \ref{prop4.1}), i.e.
 \begin{equation}\label{Ceq4.12}
 	\la^\frac{2\te}{2\te-q}\hat{D_1} < \, b^{\frac{2^*_{\mu,s}}{2^*_{\mu,s}-\te}}\left(\frac{1}{2\te} - \frac{1}{2\cdot2^*_{\mu,s}} \right)\left(S_s^H\right)^\frac{2^*_{\mu,s}\te}{2^*_{\mu,s}-\te}.
 \end{equation}
We will use the minimizers of the best constant $S_s^H$, defined in \eqref{minmimizer}, to prove the existence of solution.
%
\begin{lem}\label{lem4.2}
	For $N > 4s$ there {\color{red}exists} $\Lambda^* > 0$ such that for all $\lambda\in(0,\Lambda^*)$, 
	$$\sup\limits_{t\geq 0} J_\la(tu_\e)< c_1^* =\left( \frac{1}{2\te}- \frac{1}{2\cdot{2^*_{\mu ,s}}}\right)b^{\frac{2^*_{\mu,s}}{2^*_{\mu,s}-\te}} (S_s^H)^{\frac{2^*_{\mu,s}\te}{2^*_{\mu,s}-\te}}-\la^{\frac{2\te}{2\te-q}}\hat{D_1},$$
	with $\hat{D_1}$ as given in Proposition \ref{prop4.1}.

\proof Let $\lambda\in(0,\lambda_1)$, as defined in \eqref{Ceq4.12}, then for some positive constants $c_1$ and $c_2$
\begin{align*}
	J_\la(tu_\e) =&\; \frac{at^2}{2}\|u_\e\|^2 + \frac{bt^{2\te}}{2\theta}\|u_\e\|^{2\theta} - \frac{\lambda t^q}{q} \int\limits_{\Omega} f(x)(u_\e(x))^{q} \,dx 
	-\frac{t^{2\cdot{2^*_{\mu ,s}}}}{2\cdot{2^*_{\mu ,s}}}\|u_\e\|_0^{2\cdot{2^*_{\mu ,s}}}\\
		< &\; c_1t^2 + c_2t^{2\theta}.
\end{align*}
We can choose $t_0\in (0,1)$, such that $\sup\limits_{0\leq t\leq t_0} J_\la(tu_\e) < c_1^*$ for all $\lambda \in (0, \la_1)$. Thus it is enough to show that $\sup\limits_{t\geq t_0} J_\la(tu_\e)< c_1^*$.
\begin{align}\label{Ceq4.13o}
	\sup\limits_{t\geq t_0}J_\la(tu_\e)  &= \sup\limits_{t\geq t_0}\left( \frac{at^2}{2}\|u_\e\|^2 + \frac{bt^{2\te}}{2\theta}\|u_\e\|^{2\theta} - \frac{\lambda t^q}{q} \int\limits_{\Omega} f(x)(u_\e(x))^{q} \,dx 
	-\frac{t^{2\cdot{2^*_{\mu ,s}}}}{2\cdot{2^*_{\mu ,s}}}\|u_\e\|_0^{2\cdot{2^*_{\mu ,s}}}\right) \nonumber\\
	& \leq \sup\limits_{t\geq 0}\nu(t)-\frac{\lambda t_0^q}{q} \int\limits_{\Omega} f(x)(u_\e(x))^{q} \,dx,
\end{align}
where
\begin{equation*}
	\nu(t):=\frac{at^2}{2}\|u_\e\|^2 + \frac{bt^{2\te}}{2\theta}\|u_\e\|^{2\theta}
	-\frac{t^{2\cdot{2^*_{\mu ,s}}}}{2\cdot{2^*_{\mu ,s}}}\|u_\e\|_0^{2\cdot{2^*_{\mu ,s}}}.
\end{equation*}
Since $2< 2\te < 2\cdot2^*_{\mu,s}$, we say $\nu(0)=0$, $\nu(t)\to -\infty$ as $t\to\infty$ and $\nu(t)>0$ for small $t$.
Hence, there exists $t_\e>0$ such that $\sup\limits_{t\geq 0}\nu(t)=\nu(t_\e)$. Consequently $\nu'(t_\e)=0$, which gives
\begin{equation*}
	a\|u_\e\|^2 + bt_\e^{2\te-2}\|u_\e\|^{2\theta} = 
	t_\e^{2\cdot{2^*_{\mu ,s}}-2}\|u_\e\|_0^{2\cdot{2^*_{\mu ,s}}}.
\end{equation*}
It implies that there exists a $T_0 > 0$ such that $t_\e < T_0$. Also, we deduce that
\begin{equation*}
	a\|u_\e\|^2 \leq 
	t_\e^{2\cdot{2^*_{\mu ,s}}-2}\|u_\e\|_0^{2\cdot{2^*_{\mu ,s}}}.
\end{equation*}
Thus there exists a $T_{00} > 0$ such that $T_{00} < t_\e $.\\
Let 
\begin{equation}\label{Ceq4.113}
	\xi(t):= \frac{bt^{2\te}}{2\theta}\|u_\e\|^{2\theta}
	-\frac{t^{2\cdot{2^*_{\mu ,s}}}}{2\cdot{2^*_{\mu ,s}}}\|u_\e\|_0^{2\cdot{2^*_{\mu ,s}}}.
\end{equation}
Using $ 2\te < 2\cdot2^*_{\mu,s}$, we observe, $\xi(t)\to -\infty$ as $t\to\infty$ and $\xi(t)>0$ for small $t$.
So, there exists $t_*>0$ such that $\sup\limits_{t\geq 0}\xi(t)=\xi(t_*)$. \begin{equation}\label{Ceq4.13}
	t_*  = {\left[\frac{b\|u_\e\|^ {2\theta}}{\|u_\e\|_0^{2\cdot{2_{\mu,s}^*}}}\right]}^\frac{1}{2\cdot{2_{\mu,s}^*}-2\theta}, \quad \;
	\xi(t_*)  = \left(\frac{1}{2\theta}-\frac{1}{2\cdot 2_{\mu,s}^*}\right){\left[\frac{b\|u_\e\|^ {2\te}}{\|u_\e\|_0^{2\te}}\right]}^\frac{{2_{\mu,s}^*}}{{2_{\mu,s}^*}-\te},
\end{equation}

Also from the Proposition \ref{prop2.2} and Proposition \ref{prop2.3},
\[\|u_\e\|^ {2}\leq C^\frac{N(N-2s)}{2s(2N-\mu)}(S^H_s)^\frac{N}{2s}+ O(\e^{N-2s})\; \text{and}
\;\; \|u_\e\|_0^{2}\geq {\left[C^\frac{N}{2s}(S^H_s)^\frac{2N-\mu}{2s}- O(\e^{N})\right] }^\frac{N-2s}{2N-\mu}.\]
This implies,
\begin{align*}
	\frac{\|u_\e\|^ {2}}{\|u_\e\|_0^{2}}\leq \frac{C^\frac{N(N-2s)}{2s(2N-\mu)}(S^H_s)^\frac{N}{2s}+ O(\e^{N-2s})}{{\left[C^\frac{N}{2s}(S^H_s)^\frac{2N-\mu}{2s}- O(\e^{N})\right] }^\frac{N-2s}{2N-\mu}}
	= S^H_s \left[1+ O(\e^{N-2s})\right],
\end{align*}		 
as  $\e$ is small enough. Therefore, there exists a positive constant $c_3$, such that from \eqref{Ceq4.13}, we see
\begin{equation}\label{Ceq4.15}
	\xi(t_*)\leq\left(\frac{1}{2\theta}-\frac{1}{2\cdot 2_{\mu,s}^*}\right)b^\frac{2_{\mu,s}^*}{{2_{\mu,s}^*}-\theta}(S^H_s)^\frac{{2_{\mu,s}^*}\cdot\theta}{{2_{\mu,s}^*}-\theta} + c_3\e^{N-2s}.
\end{equation}

Let $a = \e^p$, for some $p > N-2s$. Further using \eqref{Ceq4.15} in \eqref{Ceq4.13o}, we get for some $c_4 >0$,
\begin{equation}\label{Ceq4.16}
	\begin{aligned}
		\sup\limits_{t \geq t_0}J_\la(tu_\e) & \leq\, \nu(t_\e)-\frac{\lambda t_0^q}{q} \int\limits_{\Omega} f(x)(u_\e(x))^{q} \,dx\\
		 & \leq\frac{\e^pt_\e^2}{2}\|u_\e\|^2 +\left(\frac{1}{2\theta}-\frac{1}{2\cdot 2_{\mu,s}^*}\right)b^\frac{2_{\mu,s}^*}{{2_{\mu,s}^*}-\theta}(S^H_s)^\frac{{2_{\mu,s}^*}\cdot\theta}{{2_{\mu,s}^*}-\theta} + c_3\e^{N-2s}-\frac{\lambda t_0^qm_f}{q} \|u_\e\|_{L^q(\Om)}^{q}   \\ 
		 & < \left(\frac{1}{2\theta}-\frac{1}{2\cdot 2_{\mu,s}^*}\right)b^\frac{2_{\mu,s}^*}{{2_{\mu,s}^*}-\theta}(S^H_s)^\frac{{2_{\mu,s}^*}\cdot\theta}{{2_{\mu,s}^*}-\theta} + c_4\e^{N-2s}-\frac{\lambda t_0^qm_f}{q} \|u_\e\|_{L^q(\Om)}^{q}. 
	\end{aligned}
\end{equation}
From \cite{Tuhina}, for $\e < \frac{\de}{2}$ and $N >4s$, there exists a positive constants $c_{1,s}$ such that 
\begin{equation}\label{Ceq4.17}
	\|u_\e\|_{L^q}^q \geq 
	\begin{array}{ll} 
	c_{1,s}\e^{N-\frac{q(N-2s)}{2}}, & \;\text{for}\; q > \frac{N}{N-2s}.
	\end{array} 
\end{equation}
 Using \eqref{Ceq4.16} and \eqref{Ceq4.17} with $\e = \la^{\frac{2\te}{\ba(2\te-q)}}$, where $\ba < N-2s$, for some positive constant $C_{1,s}$, we get 
 \begin{equation*}
 	\begin{aligned}
 		\sup\limits_{t \geq t_0}J_\la(tu_\e) & <\left(\frac{1}{2\theta}-\frac{1}{2\cdot 2_{\mu,s}^*}\right)b^\frac{2_{\mu,s}^*}{{2_{\mu,s}^*}-\theta}(S^H_s)^\frac{{2_{\mu,s}^*}\cdot\theta}{{2_{\mu,s}^*}-\theta} + c_4\la^{\frac{2\te}{2\te-q}}-C_{1,s}\la^{1+\left[ \frac{2\te}{\ba(2\te-q)}\right]\left[ N-\frac{q(N-2s)}{2}\right]  }.
 	\end{aligned}
 \end{equation*}
For $N > 4s$, there exists a $\ba$ 
 such that 
 \begin{equation*}
 	\sigma := -1+ \frac{2\te}{2\te-q}\left[1-\frac{1}{\ba}\left(  N-\frac{q(N-2s)}{2}\right) \right] > 0.
 \end{equation*}
 Thus, for $\la < \left( \frac{C_{1,s}}{c_4+ \hat{D_1}}\right) ^\frac{1}{\sigma}$ we have 
 \begin{equation*}
 	c_4\la^{\frac{2\te}{2\te-q}}-C_{1,s}\la^{1+\left[ \frac{2\te}{\ba(2\te-q)}\right]\left[ N-\frac{q(N-2s)}{2}\right]  } < -\la^{\frac{2\te}{2\te-q}}\hat{D_1}.
 \end{equation*}
 Define $\Lambda^* = \min\left\lbrace \la_1, \left( \frac{C_{1,s}}{c_4+ \hat{D_1}}\right) ^\frac{1}{\sigma}, \left(\frac{\de}{2}  \right)^\frac{(2\te-q)\ba}{2\te} \right\rbrace $ and $\Upsilon = {(\Lambda^*)}^{\frac{2\te}{\ba(2\te-q)}} > 0$, such that for every $\la \in (0, \Lambda^*)$ and $\e\in (0, \Upsilon)$, we have 
 	$$\sup\limits_{t \geq 0}J_\la(tu_\e) < c_1^*. \qed$$  
\end{lem}

\textbf{Proof of Theorem \ref{Cthm1.4}:} For $t$ large enough, we can write $e = tu_\e$, which along with Lemma \ref{lem4.2} gives
\begin{equation*}
	\begin{aligned}
		c_\la &\leq \max\limits_{t \in [0,1]}J_\la(te)
		< \; c_1^*,
	\end{aligned}
\end{equation*}  for all $\la \in (0, \Lambda^*)$. Thus, using Proposition \ref{prop4.1}, Lemma \ref{lem4.1} and Mountain-pass theorem \cite[Theorem 4.3.1]{badiale} we deduce that there exists a non-trivial solution of $(P_\la)$ say $u_1 \in X_0(\Om)$. Further using Lemma \ref{Clem4.4} we conclude $u_1 $ is a positive solution of $(P_\la)$. \qed

\textbf{Proof of Theorem \ref{Cthm1.3}:} The geometry of $J_\la$ for the case $q > 2\te$, implies that there exists $t_\la > 0$ such that for fixed $v \in X_0(\Om)$, we have
\begin{equation*}
	\sup\limits_{t\geq 0} J_\la(tv) = J_\la(t_\la v).
\end{equation*}
Consequently,
\begin{equation*}
	0 = \frac{at_\la}{2}\|v\|^2 + \frac{bt_\la^{2\te-1}}{2\theta}\|v\|^{2\theta} - \frac{\lambda t_\la^{q-1}}{q} \int\limits_{\Omega} f(x)(v^+(x))^{q} \,dx 
	-\frac{t_\la^{2\cdot{2^*_{\mu ,s}}-1}}{2\cdot{2^*_{\mu ,s}}}\|v^+\|_0^{2\cdot{2^*_{\mu ,s}}}.
\end{equation*}
Observe that as $\la \to \infty$, $t_\la \to 0$, and by the continuity of the functional $J_\la$, we conclude
\begin{equation*}
	\lim\limits_{\la \to \infty}\sup\limits_{t\geq 0} J_\la(tv) =\lim\limits_{\la \to \infty}J_\la(t_\la v)= 0. 
\end{equation*}
Thus, there exists $\Lambda_* > 0$ such that for all $\la \geq \Lambda_*$
\begin{equation*}
	\sup\limits_{t\geq 0} J_\la(tv) < c^*_2,
\end{equation*}
since $c^*_2 = \left( \frac{1}{2}- \frac{1}{2\cdot{2^*_{\mu ,s}}}\right) (aS_s^H)^{\frac{2N-\mu}{N-\mu+2s}} > 0$. There exists $t_0 > 0$, large enough such that $e= t_0v$, where $e$ satisfies  Lemma \ref{lem4.1} (ii). 
Now we deduce 
\begin{equation}\label{Ceq5.01} 
	\al \leq c_\la \leq \max\limits_{t \in [0,1]}J_\la(te) 
	\leq \sup\limits_{t\geq 0} J_\la(tv) < c^*_2, \;\; \text{for all} \; \la \geq \Lambda_*.
\end{equation}  
 From Proposition \ref{prop3.1}, $J_\la$ satisfies the Palais-Smale condition at the level $c_\la$ and from \eqref{Ceq5.01} we have $c_\la < c_2^*$ for all $\la \geq \Lambda_*$. Thus using Mountain-Pass theorem, there exists a non-trivial critical point of the functional $J_\la$, say $u_1 \in X_0(\Om)$ which is a non-trivial solution of $(P_\la)$. Furthermore, by Lemma \ref{Clem4.4} we conclude that $u_1 $ is a positive solution. Correspondingly using the same arguments, we can deduce similar results for the case $q = 2\te$.  \qed

\section{Multiplicity result for $\te \geq 2^*_{\mu,s}$ and $2 < q < 2^*_s$} 
In this section, we will show the existence of two positive solution of $(P_\la)$ when $2 < q < 2^*_s$ and $\te \geq 2^*_{\mu,s}$.  Here we demonstrate the proof by using the  minimization argument and Mountain-pass theorem, with condition either on $a$ or $b$ and  for $\la$ large enough. Also, here we are focusing on the case  $f > 0$.

\begin{lem}\label{Clem3.02}
	Let $2< q < 2^*_s$, then $J_\la$ is coercive for all $a >0$ and
{\small $b > \begin{cases}
		0 & \text{if}\; \te > 2^*_{\mu,s}\\
		(S_s^H)^{-2^*_{\mu,s}} & \text{if}\; \te = 2^*_{\mu,s}.
	\end{cases}
$	}
	\proof  From H\"older's inequality and Sobolev embedding we get 
	\begin{equation*}
		J_\la(u_n) \geq \frac{a}{2}\|u_n\|^2 + \frac{b}{2\theta}\|u_n\|^{2\theta} - \frac{\lambda}{q} \|f\|_{L^r}S_s^{-q/2}\|u_n\|^{q} 
		-\frac{1}{2\cdot{2^*_{\mu ,s}}}(S_s^H)^{-2^*_{\mu ,s}}\|u_n\|^{2\cdot{2^*_{\mu ,s}}}.
	\end{equation*}
	Since $2 < q < 2\cdot2^*_{\mu,s} \leq 2\te$ and by considering the assumption on $b$, we deduce that $J_\la$ is coercive.\qed
\end{lem}

 \begin{prop}\label{prop5.1}
	Let $\{u_n\}$ be a $(PS)_c$ for $J_\la$. Then $\{u_n\}$ contains a convergent subsequence for all $\la> 0$ and
	\begin{enumerate}
		\item [(i)]  $a > 0, b > (S_s^H)^{-2^*_{\mu,s}}$ when $\te = 2^*_{\mu,s}$,
		\item [(ii)] $a > 0$ and $b > \mathfrak{B}$ or $b >0$ and $a >\mathfrak{A}$ when $\te > 2^*_{\mu,s}$.
	\end{enumerate}
 where
	$\mathfrak{A}:=\frac{\te-2^*_{\mu,s}}{\te-1}\left[\frac{2^*_{\mu,s}-1}{b(\te-1)} \right]^\frac{2^*_{\mu,s}-1}{\te-2^*_{\mu,s}} (S_s^H)^{\frac{-2^*_{\mu,s}(\te-1)}{\te-2^*_{\mu,s}}}$ and $ \mathfrak{B}: = \frac{2^*_{\mu,s}-1}{\te-1}\left[\frac{\te-2^*_{\mu,s}}{a(\te-1)} \right]^\frac{\te-2^*_{\mu,s}}{2^*_{\mu,s}-1} (S_s^H)^{\frac{-2^*_{\mu,s}(\te-1)}{2^*_{\mu,s}-1}}$.
	\proof Let $\{u_n\}$ be a $(PS)_c$ for $J_\la$. Then it follows from 
Lemma \ref{Clem3.02}, $\{u_n\}$ is a bounded sequence. Thus, there exists $u \in X_0(\Om)$ such that up to a subsequence,
\begin{equation}\label{Ceq5.1}\left\{
	\begin{array}{lr}
		u_n\rightharpoonup u\; \text{in} \; X_0(\Om),\; \|u_n\|\rightarrow \ba,\;
		u_n\rightarrow u\; \text{in} \; L^{p}(\Omega)\; \text{for all}\; p\to [1,2_s^*),\\
		 \|u_n-u\|\rightarrow l ,\\
		\ds\left( \int\limits_{\Omega} \frac{(u_n^+(y))^{2^{*}_{\mu ,s}}}{|x-y|^ \mu}\,dy\right) (u_n^+)^{2^{*}_{\mu ,s}-1} \rightharpoonup \left( \int\limits_{\Omega} \frac{(u^+(y))^{2^{*}_{\mu ,s}}}{|x-y|^ \mu}\,dy\right) (u^+)^{2^{*}_{\mu ,s}-1}\; \text{weakly in} \; L^\frac{2N}{N+2s},\\
		u_n\rightarrow u \quad\text{a.e. in}\; \Om.
	\end{array}
	\right.
\end{equation}
Next, by \eqref{Ceq5.1} we have
\begin{equation}\label{Ceq5.2}
	\begin{aligned}
		o(1) =& \left\langle J'_\la(u_n)-J'_\la(u), u_n-u \right\rangle \\
		= & a\|u_n-u\|^2 + b\|u_n\|^{2\theta - 2} \|u_n-u \|^2 + b\left(\|u_n\|^{2\theta - 2} - \|u\|^{2\theta - 2}\right)  \left\langle u, u_n-u \right\rangle\\
		& -\la\int\limits_{\Om}f(x)\left( (u_n^+(x))^{q-1}-(u^+(x))^{q-1}\right) (u_n-u)(x)\,dx\\
		& - \left(\|u_n^+\|_0^{2\cdot2^*_{\mu,s}} - \|u^+\|_0^{2\cdot2^*_{\mu,s}}\right) - \|u^+\|_0^{2\cdot2^*_{\mu,s}}+\iint\limits_{\Om\times\Om}\frac{ (u^+(y))^{2^*_{\mu,s}}(u^+(x))^{2^*_{\mu,s}-1}u_n(x)}{|x-y|^\mu}\,dxdy.
	\end{aligned}
\end{equation}
By \eqref{Ceq5.1}, we have 
\begin{equation}\label{Ceq5.3}
	\lim\limits_{n \to \infty} b\left(\|u_n\|^{2\theta - 2} - \|u\|^{2\theta - 2}\right)  \left\langle u, u_n-u \right\rangle = 0.
	\end{equation}
Using Riesz representation theorem and \eqref{Ceq5.1} we have
\begin{align}\label{Ceq5.4}
	\lim\limits_{n \to \infty} \la\int\limits_{\Om}f(x)\left( (u_n^+(x))^{q-1}-(u^+(x))^{q-1}\right) (u_n-u)(x)\,dx = 0.\\ \label{Ceq5.5}
	\lim\limits_{n \to \infty} \iint\limits_{\Om\times\Om}\frac{ (u^+(y))^{2^*_{\mu,s}}(u^+(x))^{2^*_{\mu,s}-1}u_n(x)}{|x-y|^\mu}\,dxdy = \|u^+\|_0^{2\cdot2^*_{\mu,s}}.
\end{align}
Combining \eqref{Ceq5.3}, \eqref{Ceq5.4} and \eqref{Ceq5.5} in \eqref{Ceq5.2} and using Br\'ezis-Lieb lemma we get as $n \to \infty$
\begin{equation*}
	o(1) \geq a\|u_n-u\|^2 + b\|u_n\|^{2\theta - 2} \|u_n-u \|^2 - \|u_n-u\|_0^{2\cdot2^*_{\mu,s}}.
\end{equation*}
Let $v_n := u_n-u$, so by \eqref{Ceq5.1}, $\|v_n\| \to l$. Let us assume $l> 0$, otherwise we are done. Applying Br\'ezis-Lieb lemma and Sobolev embedding in the last equality we get
\begin{equation}\label{Ceq5.6}
	al^2 + bl^{2\theta } \leq (S_s^H)^{-2^*_{\mu,s}}l^{2\cdot2^*_{\mu,s}}.
\end{equation}
When $\te = 2^*_{\mu,s}$ and $b > (S_s^H)^{-2^*_{\mu,s}}$, we get from \eqref{Ceq5.6} that $l = 0$. Thus $u_n \to u$ strongly in $X_0(\Om)$.\\
When $\te > 2^*_{\mu,s}$, it follows from \eqref{Ceq5.6} and Young's inequality  
\begin{equation}\label{Ceq5.7}
	\begin{aligned}
		al^2 + bl^{2\theta } &\leq \left[ \left( \frac{a(\te-1)}{\te-2^*_{\mu,s}}\right) ^\frac{\te-2^*_{\mu,s}}{\te-1}l^\frac{2\left( \te-2^*_{\mu,s}\right) }{\te-1}\right] \left[\left( \frac{a(\te-1)}{\te-2^*_{\mu,s}}\right) ^\frac{-\left( \te-2^*_{\mu,s}\right) }{\te-1} (S_s^H)^{-2^*_{\mu,s}}l^\frac{2\te(2^*_{\mu,s}-1)}{\te-1}\right]\\
		& = al^2 + \frac{2^*_{\mu,s}-1}{\te-1}\left[\frac{\te-2^*_{\mu,s}}{a(\te-1)} \right]^\frac{\te-2^*_{\mu,s}}{2^*_{\mu,s}-1} (S_s^H)^{\frac{-2^*_{\mu,s}(\te-1)}{2^*_{\mu,s}-1}}l^{2\te},
	\end{aligned} 
\end{equation}
which contradicts the fact that $b > \mathfrak{B}$. Hence $l = 0$, i.e. $u_n \to u$ strongly in $X_0(\Om)$.
Similar analysis as in \eqref{Ceq5.7}, can be done for the case when $a > \mathfrak{A}$.\qed
\end{prop}

\textbf{Proof of Theorem \ref{Cthm1.5}:}
 First we show the existence of  least energy solution. Define 
  $$m_\la := \inf\limits_{u\in X_0(\Om)} J_\la(u).$$
  Employing  Lemma \ref{Clem3.02}, we have following deductions:
  \begin{itemize}
  	\item [(i)] For a fixed $v \in X_0(\Om)$, with $\int\limits_{\Omega} f(x)(v^+)^{q} \,dx > 0$  there exists   $\Lambda_{**} > 0$ such that for $\la > \Lambda_{**}$,  $	J_\la(v) \leq 0$. It implies $m_\la$ is well defined and $m_\la <0$.
  	\item [(ii)] We can easily choose  $\al,\rho >0$ such that $J_\la(u) > \al$ for $\|u\| \leq \rho$. 
  	\item [(iii)] There exists $0\not = e \in X_0$ such that $J_\la(e)=0$, for $\la > \Lambda_{**}$.
  \end{itemize}
By Proposition \ref{prop5.1}, $J_\la$  satisfies $(PS)_c$, say at the critical level $c= m_\la$. There exists $u_1$ in $X_0(\Om)$ such that $m_\la = \inf\limits_{u\in X_0(\Om)} J_\la(u) = J_\la(u_1). $ This implies $u_1$ is a non-trivial solution of $(P_\la)$.\\
Further, we define $$d_\la = \inf_{\xi \in \mathfrak{\Gamma}} \sup_{t \in [0,1]} J_\la(\xi(t))$$
where $\mathfrak{\Gamma}= \{\xi \in C([0,1],X_0): \xi(0)= 0,~\xi(1)= u_1\}$. Taking into account Proposition \ref{prop5.1},  and Mountain-pass theorem \cite[Theorem 2.1]{ambro}, we get a non-trivial critical point $u_2 \in X_0 $ for $J_\la$ at level $d_\la$, for $\lambda > \Lambda_{**}$.
 Solutions $u_1$ and $u_2$ are distinct, since $J_\la(u_2) = d_\la > 0 > m_\la = J_\la(u_1) $. Furthermore, by arguing as in Lemma \ref{Clem4.4}, we imply $u_1$ and $u_2$ are two positive solution of $(P_\la)$.
\qed

 \section{Multiplicity result for sublinear case $1< q \leq 2$}
  This section is devoted to the study of multiplicity of solution when we have combined effects of concave-convex non-linearties. 
   In this section, to obtain the existence of solution, we deploy the minimization arguments over a small ball in $X_0$.  Then we show the existence of second solution using the classical Mountain-pass geometry when $ 1 \leq \te < 2^*_{\mu,s}$.
\subsection{First solution}
\begin{lem}\label{lem6.1}
	There exist positive numbers $\rho $, $\lambda_0 $ and $\alpha $ such that the following holds:
	\begin{enumerate}
		\item [(i)] $J_\lambda (u)\geq  \alpha $ for any $u \in X_0$,
		with  $\|u\|  =\rho$, and for any $\lambda \in (0, \lambda_0]$.
		\item [(ii)] Let $ m_\lambda :=\inf \{J_\lambda (u) : u \in \overline{\rm B_\rho}\}$,
		where $\overline{\rm B_\rho} = \{ u \in X_0 : \lVert u \rVert \leq \rho \}$. Then ${m_\lambda < 0}$ for any ${\lambda \in (0, \lambda_0]}$.
	\end{enumerate} 
Proof (i) It follows from Case 3 of Lemma \ref{lem4.1}.\\
Proof (ii) Since $1 < q < 2 $, thus for a fixed $v\in X_0$ with $v^+\neq 0$ and $t$ small enough we have,
\begin{equation*}
	J_\lambda(tv) = \frac{at^2}{2}\|v\|^2 + \frac{bt^{2\theta}}{2\theta}{\lVert v \rVert}^{2\theta} - \frac{\lambda t^{q}}{q} \int\limits_{\Omega}f(x) (v^+)^{q} -\frac{t^{2\cdot{2^*_{\mu ,s}}}}{2\cdot{2^*_{\mu ,s} }}\iint\limits_{\Omega\times\Omega} \frac{(v^{+}(y))^{2^{*}_{\mu ,s}}(v^{+}(x))^{2^{*}_{\mu ,s}}}{|x-y|^ \mu}\,dxdy < 0.
\end{equation*}
It implies for $\|u\|$ small enough, $m_\lambda < 0$ for all $\lambda \in (0, \lambda_0]$.\qed
\end{lem}

\begin{thm}\label{Cthm6.1} There exists $\lambda_0$ such that for any $\lambda \in (0, \lambda_0]$, $(P_\la)$ has a solution $u_1 \in X_0$
	with $J_\lambda (u_1) < 0$.
	\proof Let $\rho$ and $\la_0$ be as given in Lemma \ref{lem6.1}, then we claim that there exists $u_1\in \overline{\rm B_\rho } $ such that $J_\lambda(u_1) = m_\lambda < 0$ for $\lambda \in (0, \lambda_0]$. By the definition of $m_\lambda$, there exists a minimizing sequence say $\{u_n\} \subset \overline{\rm B_\rho}$ such that
	\begin{equation}\label{Ceq6.2}
		\lim \limits_{n \to \infty} J_\lambda(u_n) = m_\lambda.
	\end{equation}
	Clearly $\{u_n\}$ is a bounded sequence in $X_0$, up to a sub-sequence,  there exists a function $u_1 \in X_0$ such that,
$u_n \rightharpoonup u_1 $ in $X_0(\Om)$, $u_n \rightharpoonup u_1$ in $ L^{2_s^*}(\Omega), u_n \rightarrow u_1$ strongly in $ L^p(\Omega)$ for all $ p\in[1,2_s^*)$, $u_n \rightarrow u_1$ a.e in $\Omega.$ We have
	\begin{equation}\label{Ceq6.3}
		\lim \limits_{n \to \infty} \int\limits_{\Omega}f(x)(u_n^+(x))^{q} \,dx =\int\limits_{\Omega}f(x)(u_1^+(x))^{q} \,dx.
	\end{equation}
	Let ${w_n = u_n-u_1}$, then by \cite[Theorem 2]{BrezLieb} and \cite[Lemma 2.2]{M.yang}, we have
	\begin{equation}\label{Ceq6.4}
		\|u_n\|^2= \|w_n\|^2+\|u_1\|^2+ o(1),\; \;
		\|u_n\|_0^{2\cdot2^{*}_{\mu ,s}} =\|w_n\|_0^{2\cdot2^{*}_{\mu ,s}} + \|u_1\|_0^{2\cdot2^{*}_{\mu ,s}} +o(1).
	\end{equation}
	Also for any $u \in B_\rho$, we deduce that the following holds
	\begin{equation}\label{Ceq6.5}
		\frac{a\|u\|^2}{2} + \frac{b\|u\|^{2\te}}{2\te}-\frac{\|u^+\|_0^{2\cdot2^*_{\mu,s}}}{2\cdot2^*_{\mu,s}} \geq 0.
	\end{equation}
	By \eqref{Ceq6.4}, we have $w_n \in \overline{\rm B_\rho}$ for $n$ sufficiently large and thus, by \eqref{Ceq6.2}, \eqref{Ceq6.3} and \eqref{Ceq6.5} it follows that, 
	\begin{align*}
		m_\lambda &= J_\lambda (u_n) + o(1)\\
		&\geq\frac{a({{\lVert w_n \rVert}^{2}+{\lVert u_1 \rVert}^{2}})}{2}+\frac{b({{\lVert w_n \rVert}^{2\te}+{\lVert u_1 \rVert}^{2\te}})}{2\theta}-\frac{\lambda }{q} \int\limits_{\Omega}f(x) (u_1^+)^{q}-\dfrac{\|w_n^{+}\|_0^{2\cdot2^{*}_{\mu ,s}}+\|u_1^{+}\|_0^{2\cdot2^{*}_{\mu ,s}}}{2\cdot{2^*_{\mu ,s} }}\\
		&\geq J_\lambda (u_1)+\frac{a\|w_n\|^2}{2} +\dfrac{b{\lVert w_n \rVert}^{2\theta}}{2\theta}-\dfrac{\|w_n^{+}\|_0^{2\cdot2^{*}_{\mu ,s}}}{2\cdot{2^*_{\mu ,s} }} +o(1)
		\geq m_\lambda.
	\end{align*}
	Noting that $\overline{\rm B_\rho }$ is a closed convex set, thus $u_1\in \overline{\rm B_\rho } $. Hence, $u_1$ is a local minimizer for $J_\lambda$, with $J_\lambda(u_1) = m_\lambda < 0$, which implies $u_1$ is a non-trivial solution.\\
	Next, it remains to prove that $u_1$ is a positive solution of $(P_\la)$. As $u_1$ is a local minimizer for $J_\lambda$, thus for any $\psi\in X_0$ and small enough $t > 0$, such that $u_1 + t\psi \in \overline{\rm B_\rho}$\; we get
	\begin{align*}
		0\leq & \,J_\lambda(u_1 + t\psi)- J_\lambda (u_1)\\
		=& \ds \frac{a\left( \|u_1+t\psi\|^2-\|u_1\|^2\right) }{2} + \frac{b\left( {\| u_1 + t\psi \|}^{2\theta}-{\lVert u_1 \rVert}^{2\theta}\right) }{2\theta}-\dfrac{\lambda }{q} \int\limits_{\Omega}f\left[ ({(u_1 + t\psi)}^+)^{q}-(u_1^+)^{q}\right]dx \\&-\dfrac{1}{2\cdot{2^*_{\mu ,s}}} \left( \|(u_1 + t\psi)^{+}\|_0^{2\cdot2^{*}_{\mu ,s}}-\|u_1^{+}\|_0^{2\cdot2^{*}_{\mu ,s}}\right) .
	\end{align*}
	So dividing by $t > 0$ and taking $t \to 0^+$, we get for any $\psi \in X_0$
	\begin{equation*}
		0 \leq \left( a+b\|u_1\|^{2\theta -2}\right) \langle u_1,\psi \rangle-\lambda \int\limits_{\Omega} f(x)(u_1^+)^{q-1}\psi dx -\iint\limits_{\Omega\times\Omega}\frac{(u_1^{+}(y))^{2^{*}_{\mu ,s}}(u_1^+(x))^{2^{*}_{\mu ,s}-1}\psi (x)}{|x-y|^ \mu}dxdy.
	\end{equation*}
	As $\psi$ was arbitrarily chosen, we infer $u_1$ is a non-trivial weak solution.
Employing Lemma \ref{Clem4.4}, we imply $u_1$ is a positive solution of $(P_\la)$.\qed
\end{thm}	 
Thus from Theorem \ref{Cthm6.1}, we conclude Remark \ref{Crem1.1}.

 \subsection{Second solution}
 	
 We obtain the second solution with the help of minimizers of the best constant $S_s^H$, defined in \eqref{minmimizer}.\\
Choose $\la_1 > 0$, in such manner that for $\la \in (0, \la_1)$, we have $c^*_3 >0$ (as defined in Proposition \ref{P.S}). Let 
 \begin{equation}\label{Ceq2.18n}
 	\la_* = \min \{\la_0, \la_1\}.
 \end{equation}
 Next we recall $u_1$, which is the local minimum $J_\la$ and the minimizer $u_\e$ from \eqref{Ceq2.4}. Also we recollect the definition of \,$\Xi$ from Section 3 and the existence of $\de > 0$ such that $B_{4\de}(0) \subset \Xi \subseteq \Om$ and $f(x) > 0$ for all $x \in B_{2\de}(0)$. It implies that there exists a $m_f > 0$, such that $f(x) > m_f$ for all $x \in B_{2\de}(0)$. The proof of the following Lemma is similar to that of \cite[Lemma 4.2]{Div}. Thus we just give a sketch of the proof for the readers' convenience. 
 \begin{lem}\label{lem6.3}
 	Let $\mu < \min\{4s, N\}$, then there exists $\Lambda^{**} > 0$ such that for all $0 <\la < \Lambda^{**}$ we obtain
 	\begin{equation*}
 		\sup\limits_{t \geq 0}J_\la(u_1+tu_\e) < c^*_3.
 	\end{equation*}
 	\proof Consider
 	\begin{equation}\label{Ceq2.19}
 		J_\la(u_1+tu_\e) = \frac{a\|u_1+tu_\e\|^2}{2}+ \frac{b\|u_1+tu_\e\|^{2\te}}{2\te}-\frac{\la}{q}\int\limits_{\Om}f(x)(u_1+tu_\e)^q\;dx-\frac{\|u_1+tu_\e\|_0^{2\cdot2^*_{\mu,s}}}{2\cdot2^*_{\mu,s}}.
 	\end{equation}
 	For some $\ba \in [0,2\pi)$ we have
 	\begin{align}\label{Ceq2.20}
 		\|u_1+tu_\e\|^2 &= \|u_1\|^2 + \|tu_\e\|^2 + 2\left\langle u_1, tu_\e \right\rangle \nonumber\\
 		& = \|u_1\|^2 + \|tu_\e\|^2 + 2\|u_1\|\|tu_\e\|\cos\beta,
 	\end{align}
 	which implies $\|u_1+tu_\e\|^{2\te } = \left( \|u_1\|^2 + \|tu_\e\|^2 + 2\|u_1\|\|tu_\e\|\cos\ba\right) ^{\te}.$\\
 	Additionally for all $y \geq 0, \ba \in [0, 2\pi)$, there exists a uniform $R > 0$ such that following inequality holds
 	\begin{equation*}
 		\left( 1+ y^2 + 2y \cos\ba\right) ^\te \leq \begin{cases}
 			1+ y^{2\te} + 2\te y \cos\ba +Ry^2 & \te \in [1,\frac{3}{2}),\\
 			1+ y^{2\te} + 2\te y \cos\ba +R(y^{2\te -1} + y^2) & \te \geq \frac{3}{2}.
 		\end{cases}
 	\end{equation*}
 	Putting $y = \frac{\|tu_\e\|}{\|u_1\|}$ in the above inequality, we have the following estimate for some $c_0 > 0$
 	\begin{equation}
 		\begin{aligned}\label{Ceq2.21}
 			\|u_1+tu_\e\|^{2\te }
 			&\leq \|u_1\|^{2\te} + \|tu_\e\|^{2\te} + 2\te\|u_1\|^{2\te-1}\|tu_\e\|\cos\ba + Rc_0(\|tu_\e\|^{2\te-1} + \|tu_\e\|^{2})\\
 			& = \|u_1\|^{2\te} + \|tu_\e\|^{2\te} + 2\te\|u_1\|^{2\te-2}\left\langle u_1, tu_\e\right\rangle  + Rc_0(\|tu_\e\|^{2\te-1} + \|tu_\e\|^{2}).
 		\end{aligned}
 	\end{equation}
 	From \cite[Lemma 4.2]{Div}, we get the estimate of $\|u_1+ tu_\e\|_0^{2\cdot2^*_{\mu,s}}$ for the case when $2^*_{\mu,s} > 3$ and $2 < 2^*_{\mu,s} \leq 3$, as follows for all $\tau\in(0,1)$ and for some $\hat{C} > 0$,
 	\begin{equation}\label{Ceq2.22}
 		\begin{aligned}
 			\|u_1+ tu_\e\|_0^{2\cdot2^*_{\mu,s}} \geq& \|u_1\|_0^{2\cdot2^*_{\mu,s}} + \|tu_\e\|_0^{2\cdot2^*_{\mu,s}}	\\ &+ 2\cdot2^*_{\mu,s}t^{2\cdot2^*_{\mu,s}-1}\hat{C}\iint\limits_{\Om\times\Om}\frac{(u_\e(x))^{2^*_{\mu,s}}(u_\e(y))^{2^*_{\mu,s}-1}u_1(y)}{|x-y|^\mu} \,dxdy\\
 			&+ 2\cdot2^*_{\mu,s}t\iint\limits_{\Om\times\Om}\frac{(u_1(x))^{2^*_{\mu,s}}(u_1(y))^{2^*_{\mu,s}-1}u_\e(y)}{|x-y|^\mu}\,dxdy - O(\e^{\left( \frac{2N-\mu}{4}\right)\tau }).
 		\end{aligned}
 	\end{equation}
 	Putting together \eqref{Ceq2.20}, \eqref{Ceq2.21} and \eqref{Ceq2.22} in \eqref{Ceq2.19}, we get
 	\begin{equation}\label{Ceq2.23}
 		\begin{aligned}
 			J_\la(u_1+tu_\e) \leq& \frac{a\|u_1\|^2}{2}+ \frac{a\|tu_\e\|^2}{2}+ a\left\langle u_1, tu_\e\right\rangle + \frac{b\|u_1\|^{2\te}}{2\te}+ \frac{b\|tu_\e\|^{2\te}}{2\te}+ b\|u_1\|^{2\te-2}\left\langle u_1, tu_\e\right\rangle \\
 			&+ \frac{bRc_0}{2\te}(\|tu_\e\|^{2\te-1} + \|tu_\e\|^{2})-\frac{\la}{q}\int\limits_{\Om}f(u_1+tu_\e)^q\;dx-\frac{\|u_1\|_0^{2\cdot2^*_{\mu,s}}}{2\cdot2^*_{\mu,s}} - \frac{\|tu_\e\|_0^{2\cdot2^*_{\mu,s}}}{2\cdot2^*_{\mu,s}} \\
 			&-\hat{C} t^{2\cdot2^*_{\mu,s}-1}\iint\limits_{\Om\times\Om}\frac{(u_\e(x))^{2^*_{\mu,s}}(u_\e(y))^{2^*_{\mu,s}-1}u_1(y)}{|x-y|^\mu}\,dxdy \\
 			&- t\iint\limits_{\Om\times\Om}\frac{(u_1(x))^{2^*_{\mu,s}}(u_1(y))^{2^*_{\mu,s}-1}u_\e(y)}{|x-y|^\mu}\,dxdy + O(\e^{\left( \frac{2N-\mu}{4}\right)\tau }).
 		\end{aligned}
 	\end{equation}
 	Taking $\tau =\frac{2}{2^*_{\mu,s}}$, and using the fact that $u_1$ solves $(P_\la)$, i.e. $\left\langle J'_\la(u_1), tu_\e\right\rangle  = 0$ and $J_\la(u_1) < 0$, in \eqref{Ceq2.23} we get
 	\begin{equation}\label{Ceq2.24}
 		\begin{aligned}
 			J_\la(u_1+tu_\e) <& \frac{a\|tu_\e\|^2}{2}+  \frac{b\|tu_\e\|^{2\te}}{2\te} + \frac{bRc_0}{2\te}(\|tu_\e\|^{2\te-1} + \|tu_\e\|^{2})\\
 			&-\la\int\limits_{\Om}f(x)\left(\int\limits_0^{tu_\e} (u_1+s)^{q-1}-u_1^{q-1}\,ds\right) \;dx  - \frac{\|tu_\e\|_0^{2\cdot2^*_{\mu,s}}}{2\cdot2^*_{\mu,s}} \\
 			&-\hat{C} t^{2\cdot2^*_{\mu,s}-1}\iint\limits_{\Om\times\Om}\frac{(u_\e(x))^{2^*_{\mu,s}}(u_\e(y))^{2^*_{\mu,s}-1}u_1(y)}{|x-y|^\mu}\,dxdy + o(\e^{ \frac{N-2s}{2} }). 
 		\end{aligned}
 	\end{equation}
 	Also one can show there exists $r_1 > 0$ such that
 	\begin{equation}\label{Ceq2.25}
 		\iint\limits_{\Om\times\Om}\frac{(u_\e(x))^{2^*_{\mu,s}}(u_\e(y))^{2^*_{\mu,s}-1}u_1(y)}{|x-y|^\mu}\,dxdy \geq r_1\e^{\frac{N-2s}{2}}.
 	\end{equation}
 	Using \eqref{Ceq2.25} in \eqref{Ceq2.24}, we get
 	\begin{equation}\label{Ceq2.26}
 		\begin{aligned}
 			J_\la(u_1+tu_\e) < & \frac{a\|tu_\e\|^2}{2}+  \frac{b\|tu_\e\|^{2\te}}{2\te} + \frac{bRc_0}{2\te}(\|tu_\e\|^{2\te-1} + \|tu_\e\|^{2})\\
 			& - \frac{\|tu_\e\|_0^{2\cdot2^*_{\mu,s}}}{2\cdot2^*_{\mu,s}} 
 			-\hat{C} t^{2\cdot2^*_{\mu,s}-1}r_1\e^{\frac{N-2s }{2}}+ o(\e^{ \frac{N-2s}{2} }). 
 		\end{aligned}
 	\end{equation}
 	Define
 	\begin{equation*}
 		\mathcal{A}(t):=\frac{a\|tu_\e\|^2}{2}+  \frac{b\|tu_\e\|^{2\te}}{2\te} + \frac{bRc_0}{2\te}(\|tu_\e\|^{2\te-1} + \|tu_\e\|^{2})\\
 		- \frac{\|tu_\e\|_0^{2\cdot2^*_{\mu,s}}}{2\cdot2^*_{\mu,s}} 
 		-\hat{C} t^{2\cdot2^*_{\mu,s}-1}r_1\e^{\frac{N-2s }{2}}.
 	\end{equation*}
 	From the geometry of $\mathcal{A}$, there exists $t_\e>0$ such that $\sup\limits_{t\geq 0}\mathcal{A}(t)=\mathcal{A}(t_\e)$. Consequently $\mathcal{A}'(t_\e)=0$, which implies
 	that there exist $T_0, T_{00} > 0$, such that $T_{00} < t_\e < T_0$. 
 	Let us define
 	\begin{equation*}
 		\xi(t):= \frac{at^{2}}{2}\|u_\e\|^{2}
 		-\frac{t^{2\cdot{2^*_{\mu ,s}}}}{2\cdot{2^*_{\mu ,s}}}\|u_\e\|_0^{2\cdot{2^*_{\mu ,s}}}.
 	\end{equation*}
 	Thus, there exists $0 < t_*:= {\left[\frac{a\|u_\e\|^ {2}}{\|u_\e\|_0^{2\cdot{2_{\mu,s}^*}}}\right]}^\frac{1}{2\cdot{2_{\mu,s}^*}-2}$ such that $\sup\limits_{t\geq 0}\xi(t)=\xi(t_*)$. 
 	From Proposition \ref{prop2.2} and Proposition \ref{prop2.3},
 	\begin{equation}\label{Ceq2.27}
 		\xi(t_*)\leq\left(\frac{1}{2}-\frac{1}{2\cdot 2_{\mu,s}^*}\right)(aS^H_s)^\frac{2_{\mu,s}^*}{{2_{\mu,s}^*}-1} + O(\e^{N-2s}),
 	\end{equation}
 	as $\e$ is small enough.
 	Thus, from \eqref{Ceq2.26} and \eqref{Ceq2.27}, we get
 	\begin{equation}\label{Ceq7.15}
 		\begin{aligned}
 			J_\la(u_1+ tu_\e) \leq & \left(\frac{1}{2}-\frac{1}{2\cdot 2_{\mu,s}^*}\right)(aS^H_s)^\frac{2_{\mu,s}^*}{{2_{\mu,s}^*}-1} + O(\e^{N-2s}) + \frac{b\|T_0u_\e\|^{2\te}}{2\te} + \frac{bRc_0}{2\te}\|T_0u_\e\|^{2\te-1} \\
 			&+ \frac{bRc_0}{2\te}\|T_0u_\e\|^{2}
 			-\hat{C} T_{00}^{2\cdot2^*_{\mu,s}-1}r_1\e^{\frac{N-2s }{2}} + o(\e^\frac{N-2s}{2}).
 		\end{aligned}
 	\end{equation}
 	Let $b = \e^p$, for some $p > N-2s$. Consider
 	\begin{equation*}
 		O(\e^{N-2s}) + \frac{\e^p\|T_0u_\e\|^{2\te}}{2\te} + \frac{\e^pRc_0}{2\te}\|T_0u_\e\|^{2\te-1} + \frac{\e^pRc_0}{2\te}\|T_0u_\e\|^{2}
 		+ o(\e^\frac{N-2s}{2}) \leq c_1\e^\ba
 	\end{equation*}
 	where $\ba > \frac{N-2s}{2}$ and $c_1$ is a positive constant. 
 	 Assuming $\e = \la^{\frac{2}{\ba(2-q)}}$ and using the above inequality in \eqref{Ceq7.15}, we get 
 	\begin{equation*}
 		\begin{aligned}
 			J_\la(u_1+ tu_\e) \leq  \left(\frac{1}{2}-\frac{1}{2\cdot 2_{\mu,s}^*}\right)(aS^H_s)^\frac{2_{\mu,s}^*}{{2_{\mu,s}^*}-1} +c_1 \la^{\frac{2}{2-q}} -\hat{C} T_{00}^{2\cdot2^*_{\mu,s}-1}r_1\left(\la^{\frac{2}{2-q}} \right)^ \frac{N-2s }{2\ba}.
 		\end{aligned}
 	\end{equation*}
 	Thus, for $\la < \left( \frac{\hat{C} T_{00}^{2\cdot2^*_{\mu,s}-1}r_1}{c_1+ \hat{D}}\right) ^\frac{(2-q)\ba}{2\ba-(N-2s)}:= \Lambda$ we have 
 	\begin{equation*}
 		c_1 \la^{\frac{2}{2-q}} -\hat{C} T_{00}^{2\cdot2^*_{\mu,s}-1}r_1\left(\la^{\frac{2}{2-q}} \right)^ \frac{N-2s }{2\ba} < -\la^{\frac{2}{2-q}}\hat{D}.
 	\end{equation*}
 	Define $\Lambda^{**} = \min\left\lbrace \la_*, \Lambda\right\rbrace $ and $\e_* = (\Lambda^{**})^{\frac{2}{\ba(2-q)}} > 0$, such that for every $\la \in (0, \Lambda^{**})$ and $\e\in (0, \e_*)$, we have 
 	\begin{equation*}
 		\sup\limits_{t \geq 0}J_\la(u_1+tu_\e) < c_3^*.
 	\end{equation*}\qed
 	
  \end{lem}

 \begin{lem}\label{lem3.10}
 	Let $\mu \geq 4s$ and $\frac{N}{N-2s} \leq q <2$, then there {\color{red}exists} a $\tilde{\Lambda}^{**} > 0$ such that for all $0 <\la < \tilde{\Lambda}^{**}$ we obtain
 	\begin{equation*}
 		\sup\limits_{t \geq 0}J_\la(tu_\e) < c^*_3.
 	\end{equation*}
 	\proof Let $\lambda\in(0,\la_*)$, as defined in \eqref{Ceq2.18n}, then following the same argument as in Lemma \ref{lem4.2}, we get 
 	\begin{align}\label{Ceq2.28}
 		\sup\limits_{t\geq t_0}J_\la(tu_\e)  &= \sup\limits_{t\geq t_0}\left( \frac{at^2}{2}\|u_\e\|^2 + \frac{bt^{2\te}}{2\theta}\|u_\e\|^{2\theta} - \frac{\lambda t^q}{q} \int\limits_{\Omega} f(x)(u_\e(x))^{q}dx 
 		-\frac{t^{2\cdot{2^*_{\mu ,s}}}}{2\cdot{2^*_{\mu ,s}}}\|u_\e\|_0^{2\cdot{2^*_{\mu ,s}}}\right) \nonumber\\
 		& \leq \sup\limits_{t\geq 0}\nu(t)-\frac{\lambda t_0^q}{q} \int\limits_{\Omega} f(x)(u_\e(x))^{q}dx.
 	\end{align}
 	Proceeding as in Lemma \ref{lem4.2} up to \eqref{Ceq4.113}, we conclude there exists $t_\e$ such that $\sup\limits_{t\geq 0}\nu(t) = \nu(t_\e)$, also we define
 	\begin{equation*}
 		\mathcal{K}(t):= \frac{at^{2}}{2}\|u_\e\|^{2}
 		-\frac{t^{2\cdot{2^*_{\mu ,s}}}}{2\cdot{2^*_{\mu ,s}}}\|u_\e\|_0^{2\cdot{2^*_{\mu ,s}}}.
 	\end{equation*}
 	Therefore, there exists $0 < \tilde{t}_* := {\left[\frac{a\|u_\e\|^ {2}}{\|u_\e\|_0^{2\cdot{2_{\mu,s}^*}}}\right]}^\frac{1}{2\cdot{2_{\mu,s}^*}-2}$ such that $\sup\limits_{t\geq 0}\mathcal{K}(t)=\mathcal{K}(\tilde{t}_*)$. 
 	Using Proposition \ref{prop2.2} and Proposition \ref{prop2.3},
 	we conclude there exists a positive constant $c$, such that
 	\begin{equation}\label{Ceq2.30}
 		\mathcal{K}{\color{red}(\tilde{t}_*)}\leq\left(\frac{1}{2}-\frac{1}{2\cdot 2_{\mu,s}^*}\right)(aS^H_s)^\frac{{2_{\mu,s}^*}}{{2_{\mu,s}^*}-1} + c\e^{N-2s}.
 	\end{equation}
 	In \eqref{Ceq2.28}, assuming $b = \e^p$ where $p > N-2s$ and using \eqref{Ceq2.30}, we get for some $C>0$
 	\begin{equation}\label{Ceq2.31}
 		\begin{aligned}
 			\sup\limits_{t \geq t_0}J_\la(tu_\e) 
 			& \leq \left(\frac{1}{2}-\frac{1}{2\cdot 2_{\mu,s}^*}\right)(aS^H_s)^\frac{{2_{\mu,s}^*}}{{2_{\mu,s}^*}-1} + C\e^{N-2s} -\frac{\la t_0^qm_f}{q} \int\limits_{\Omega} (u_\e(x))^{q}dx.
 		\end{aligned}
 	\end{equation}
 	For $\e < \frac{\de}{2}$, there exist positive constants $c_{1,s}, c_{2,s}$ such that
 	\begin{equation}\label{Ceq2.32}
 		\|u_\e\|_{L^q}^q \geq \left\{
 		\begin{array}{ll} c_{1,s}\e^{N-\frac{q(N-2s)}{2}} & q > \frac{N}{N-2s},\\
 			c_{2,s}\e^\frac{N}{2}|\log\e| & q = \frac{N}{N-2s}.\\
 		\end{array} 
 		\right. 
 	\end{equation}
 	Using \eqref{Ceq2.31} and \eqref{Ceq2.32} with $\e = \left( \la^{\frac{2}{2-q}}\right)^\frac{1}{N-2s}$, for some constants $C_{1,s}, C_{2,s}>0$, we get \\
 	\textbf{Case (i):} When $q >\frac{N}{N-2s}$ 
 	\begin{equation*}
 		\begin{aligned}
 			\sup\limits_{t \geq t_0}J_\la(tu_\e) & \leq \left(\frac{1}{2}-\frac{1}{2\cdot 2_{\mu,s}^*}\right)(aS^H_s)^\frac{{2_{\mu,s}^*}}{{2_{\mu,s}^*}-1} + C\la^{\frac{2}{2-q}}-C_{1,s}\la^{1+\left[ \frac{2}{(2-q)(N-2s)}\right]\left[ N-\frac{q(N-2s)}{2}\right]  }.
 		\end{aligned}
 	\end{equation*}
 	Note that $q > \frac{N}{N-2s}$, if and only if 
 	$\sigma := -1+ \frac{2}{2-q}\left[1-\frac{1}{N-2s}\left(  N-\frac{q(N-2s)}{2}\right) \right] > 0.$	
 	Therefore, there exists $\upsilon_1 := \left( \frac{C_{1,s}}{C+ \hat{D}}\right) ^\frac{1}{\sigma}$, such that for $\la < \upsilon_1$ we get 
 		$\sup\limits_{t \geq t_0}J_\la(tu_\e) < c_3^*.$\\
 	\textbf{Case (ii):} When $q =\frac{N}{N-2s}$\\
 	As $\la \to 0$, then $|\log\left( \la^{\frac{2}{2-q}}\right)^\frac{1}{N-2s}| \to \infty$, thus, there exists $\upsilon_2$, such that for $\la \in (0, \upsilon_2)$ we have
 	\begin{equation*}
 		C\la^{\frac{2}{2-q}}-C_{2,s} \la^{1+ \frac{N}{(2-q)(N-2s)}}|\log\left( \la^{\frac{2}{2-q}}\right)^\frac{1}{N-2s}| < -\la^{\frac{2}{2-q}}\hat{D}.
 	\end{equation*}
 	Define $\tilde{\Lambda}^{**} = \min\left\lbrace \la_*, \upsilon_1, \upsilon_2,  \left(\frac{\de}{2}  \right)^\frac{(2-q)(N-2s)}{2} \right\rbrace $ and $\e_{**} = (\tilde{\Lambda}^{**})^{\frac{2}{(2-q)(N-2s)}} > 0$, such that for every $\la \in (0, \tilde{\Lambda}^{**})$ and $\e\in (0, \e_{**})$, we have 
 	\begin{equation*}
 		\sup\limits_{t \geq 0}J_\la(tu_\e) < c_3^*. 
 	\end{equation*}
  \qed
\end{lem}
 
\textbf{Proof of Theorem \ref{Cthm1.5} (i):} Let $\la \in (0, \Lambda^{**})$ and $\e \in(0, \e_*)$, then from Lemma \ref{lem4.1}, $J_\la$ satisfies the geometry of the Mountain-pass lemma, thus there exists a sequence $\{u_n\}$ such that
 $J_\la(u_n) \to c_\la$ and $J'_\la(u_n) \to 0$, where
 \begin{equation*}
 	c_\la := \inf\limits_{h \in \Gamma}\max\limits_{t \in [0,1]}J_\la(h(t)),\;\; \text{where} \;
 	\Gamma = \{h \in C([0,1], X_0(\Om)): h(0) = u_1\;\text{and}\; h(1) = u_1+t_0u_\e\}.
 \end{equation*}
 Consequently by Proposition \ref{P.S} and Lemma \ref{lem6.3} we get
 $$0 < \rho < c_\la \leq \max\limits_{t \in [0,1]} J_\la(u_1 + tt_0u_\e ) \leq \sup\limits_{t \geq 0} J_\la(u_1 + tu_\e ) < c_3^*.$$
 Thus there exists $u_2 \in X_0(\Om)$ such that $u_n \to u_2$ in $X_0(\Om)$ i.e. $u_2$ is the non-trivial critical point of $J_\la$. The solutions obtained are distinct, since $J_\la(u_2) = c_\la > 0 > m_\la = J_\la(u_1) $.
 
 \textbf{Proof of Theorem \ref{Cthm1.5} (ii):} Similarly by the same argument for $\la \in (0, \tilde{\Lambda}^{**})$ and $\e \in(0, \e_{**})$, we set the minimax level
  \begin{equation*}
 	\tilde{c}_\la := \inf\limits_{g \in \Gamma}\max\limits_{t \in [0,1]}J_\la(g(t)),\, \text{where}\;
 	\Gamma = \{g \in C([0,1], X_0(\Om)): g(0) = 0\;\text{and}\; g(1) = t_0u_\e\}.
 \end{equation*}
 Therefore by Proposition \ref{P.S} and Lemma \ref{lem3.10} we get
 $ \rho < \tilde{c}_\la 
 \leq \sup\limits_{t \geq 0} J_\la(tu_\e ) < c_3^*.$
 Thus there exists $\tilde{u}_2 \in X_0(\Om)$ such that $u_n \to \tilde{u}_2$ in $X_0(\Om)$. Note that, the solutions are distinct, as $J_\la(\tilde{u}_2) = c_\la' > 0 > m_\la = J_\la(u_1) $.\\
 Next reasoning as in Lemma \ref{Clem4.4}, we get the desired result. \qed

 \subsection{Case when q =2}
 We will be studying the problem $(P_\la)$ with $q= 2$. Thus, the energy functional is defined as,
 \begin{equation*}
 	\begin{aligned}
 		J_\lambda (u)= &\frac{a}{2}\|u\|^2 + \frac{b}{2\theta}\|u\|^{2\theta} - \frac{\lambda}{2} \int\limits_{\Omega} f(x)(u^+(x))^{2} \,dx 
 		-\frac{1}{2\cdot{2^*_{\mu ,s}}}\|u^+\|_0^{2\cdot{2^*_{\mu ,s} }}.
 	\end{aligned}
 \end{equation*}
 \begin{lem}\label{lem3.12}
 	Let $N \geq 4s$, then there exists $\Upsilon > 0$ such that for $\e \in (0, \Upsilon)$ we have
 	\begin{equation*}
 		\sup\limits_{t \geq 0}J_\la(tu_\e) < c_4^* := \left( \frac{1}{2}- \frac{1}{2\cdot2^*_{\mu,s}}\right) (aS_s^H)^{\frac{2N-\mu}{N-\mu+2s}}.
 	\end{equation*}
 \proof	Proceeding as in Lemma \ref{lem3.10} up to \eqref{Ceq2.31} we have for all $\la > 0$ and for some $C >0$
 	\begin{equation}\label{Ceq3.26}
 		\begin{aligned}
 			\sup\limits_{t \geq t_0}J_\la(tu_\e) &  \leq \left(\frac{1}{2}-\frac{1}{2\cdot 2_{\mu,s}^*}\right)(aS^H_s)^\frac{{2_{\mu,s}^*}}{{2_{\mu,s}^*}-1} + C\e^{N-2s} -\frac{\la t_0^qm_f}{2} \int\limits_{\Omega} (u_\e)^{2}dx.
 		\end{aligned}
 	\end{equation}
 	Using Proposition \ref{prop2.2}, we have for some positive constant $C_s$
 	\begin{equation}\label{Ceq3.27}
 		\|u_\epsilon\|_{L^2}^2 \geq \left\{
 		\begin{array}{ll} C_s\epsilon^{2s}+ O(\epsilon^{N-2s}) & N > 4s,\\
 			C_s\epsilon^{2s}|\log(\epsilon)|+ O(\epsilon^{2s}) & N= 4s.\\
 		\end{array} 
 		\right. 
 	\end{equation}
 	Further using \eqref{Ceq3.26} and \eqref{Ceq3.27} we get\\ 
 	\textbf{Case (i)} When $N > 4s$
 	\begin{equation*}
 		\begin{aligned}
 			\sup\limits_{t \geq t_0}J_\la(tu_\e) &  \leq \left(\frac{1}{2}-\frac{1}{2\cdot 2_{\mu,s}^*}\right)(aS^H_s)^\frac{{2_{\mu,s}^*}}{{2_{\mu,s}^*}-1} + C\e^{N-2s} -\la C_s\epsilon^{2s}.
 		\end{aligned}
 	\end{equation*}
 	Thus we can choose $\Upsilon_1 > 0$ such that for $\e \in (0, \Upsilon_1)$, $C\e^{N-2s} -\la C_s\epsilon^{2s} < 0$.\\
 	\textbf{Case (ii)} When $N = 4s$
 	\begin{equation*}
 		\begin{aligned}
 			\sup\limits_{t \geq t_0}J_\la(tu_\e) &  \leq \left(\frac{1}{2}-\frac{1}{2\cdot 2_{\mu,s}^*}\right)(aS^H_s)^\frac{{2_{\mu,s}^*}}{{2_{\mu,s}^*}-1} + C\e^{N-2s} -\la C_s\epsilon^{2s}|\log(\epsilon)|.
 		\end{aligned}
 	\end{equation*}
 	As $\e \to 0,$ then $|\log(\epsilon)| \to \infty$. Thus we can choose $\Upsilon_2 > 0$ such that for $\e \in (0, \Upsilon_2)$, $C\e^{N-2s} -\la C_s\epsilon^{2s}|\log(\epsilon)| < 0$. Define $\Upsilon = \min\{\Upsilon_1, \Upsilon_2\}$, hence for all $\la > 0$ and $\e \in (0, \Upsilon)$ we have $\sup\limits_{t \geq 0}J_\la(tu_\e) < c_4^*.$\qed

 \end{lem}
\textbf{Proof of Theorem \ref{Cthm1.5} (iii):}
 Thus by Lemma \ref{lem4.1} and \ref{lem3.12} we conclude about the existence of non-trivial solution of $(P_\la)$, say $u_1 \in X_0(\Om)$. By Lemma \ref{Clem4.4}, we imply $u_1 >0$. \qed

 \noindent
 \textbf{Acknowledgement.}
 The research of  second author is  supported by a grant of UGC (India)  with  JRF grant number: June 18-414344.

\end{document}